\documentclass[11pt]{article}

\usepackage{graphics}

\usepackage{latexsym}
\usepackage{amsmath}
\usepackage{amssymb}
\usepackage[dvips]{epsfig}

\newtheorem{thm}{Theorem}[section]
\newtheorem{cor}[thm]{Corollary}
\newtheorem{lem}[thm]{Lemma}
\newtheorem{prop}[thm]{Proposition}
\newtheorem{defn}[thm]{Definition}
\newtheorem{rem}[thm]{Remark}

 \numberwithin{equation}{section}

\newcommand{\bt}{\begin{theorem}}
\newcommand{\bl}{\begin{lemma}}
\newcommand{\el}{\end{lemma}}
\newcommand{\et}{\end{theorem}}

\newcommand{\bn}{\begin{eqnarray}}
\newcommand{\en}{\end{eqnarray}}
\newcommand{\bnn}{\begin{eqnarray*}}
\newcommand{\enn}{\end{eqnarray*}}

\makeatletter      
\@addtoreset{equation}{section}
\makeatother       

\title{On 3D Lagrangian Navier-Stokes $\alpha$ model with a Class of Vorticity-Slip Boundary conditions
\footnote{This research is
supported in part by NSFC 10971174, and Zheng Ge Ru Foundation,
and Hong Kong RGC Earmarked Research Grants CUHK-4041/11P,
CUHK-4042/08P and a Focus Area Grant from The Chinese University of Hong Kong.}}

\author{Yuelong Xiao \\
{\small School of Mathematics and Computational Science}\\
{\small Xiangtan University}\\
{\small Xiangtan 411105, P.R.China}\\
{\small The Institute of Mathematical Sciences}\\
{\small The Chinese University of Hong Kong}\\
{\small Shatin, N.T., Hong Kong}\\
{xyl@xtu.edu.cn} \\ \\
{Zhouping Xin} \\
{\small The Institute of Mathematical Sciences and Department of Mathematics}\\
{\small The Chinese University of Hong Kong}\\
{\small Shatin, N.T., Hong Kong}\\
{zpxin@ims.cuhk.edu.hk}}

\date{ }
\begin{document}
\maketitle

\vskip 1cm

\[Abstract\]
\\
This paper concerns the 3-dimensional Lagrangian Navier-Stokes $\alpha$ model and
the limiting Navier-Stokes system on smooth bounded domains with a class of vorticity-slip
boundary conditions and the Navier-slip boundary conditions. It establishes the spectrum
properties and regularity estimates of the associated Stokes operators, the local
well-posedness of the strong solution and global existence of weak solutions for initial
boundary value problems for such systems. Furthermore, the vanishing $\alpha$ limit to
a weak solution of the corresponding initial-boundary value problem of the Navier-Stokes system
is proved and a rate of convergence is shown for the strong solution.

\section{Introduction}

The Lagragian Navier-Stokes $\alpha$ model (LNS-$\alpha$) as
a regularization system of the Navier-Stokes equations (NS) is given by
\begin{eqnarray}
&& \partial_t v - \Delta v + T_\alpha v\cdot\nabla v  + \nabla
(T_\alpha v)^T\cdot v +
\nabla p=0\\
&& \nabla\cdot v = 0
\end{eqnarray}
which describes large scale fluid motions in the turbulence theory,
where $T_\alpha v = u$ is a filtered version of the velocity $v$ determined usually by
\begin{eqnarray}
&& u - \alpha \Delta u = v \\
&& \nabla\cdot u = 0
\end{eqnarray}
with $\alpha > 0$ being a constant. This filter $u$ is also called the
averaged velocity. The system can be regarded as a system
for this filter, and is also called the
Lagrangian averaged Navier-Stokes equations (LANS). The ideal case, called
the Lagrangian averaged Euler equations
(LAE) or Camass-Holm equations, was first introduced in \cite{CF8,HM1}.
The viscosity was added in \cite{CF1,CF2,H9}
yielding the LANS which is sometime called viscous Camass-Holm equations.
\\

The global well-posedness for the LANS was first obtained in
\cite{FHT} for periodic boundary conditions. The convergence
of its solutions to that of the NS equations and the
continuity of attractors when $\alpha\rightarrow 0$ are also considered there.
\\

For bounded domains, the situation becomes more complicated since
the LANS is a 4th odder system for the filter $u$, and only the
no-slip boundary condition $u = 0$ on the boundary was considered
by \cite{MS} under the assumption that $Au = -P\Delta u = 0$ on the
boundary with $P$ being the Leray projection operator. The boundary effects related to such a boundary
condition were analyzed in \cite{H03}. We also refer
\cite{FHT,GH6,H03,MS}
for more details along this line.
\\

On the other hand, the LNS-$\alpha$ model emphasizes the system (1.1)-(1.4) as equations for
the physical velocity $v$, which is a regularized system of the NS equations by
filtering some part of the nonlinearity through a global quantity which is then
called filtered velocity (see \cite{GH6} and the references therein). There are many
filtered formulations, which thus lead to many $\alpha$ models(see \cite{CT,GH8} for
instance). It is also mentioned in \cite{C08} in the stochastic Lagrangian derivation
of (1.1), (1.2) that any translation-invariant filter $u = T_\alpha v$ may be adaptable.
\\

Although, there is no any serious difference between the two aspects
for the equations (1.1), (1.2) filtered by (1.3), (1.4) in domains without boundary,
the situation may be different for domains with boundaries.
To our knowledge, very little is known to the LNS-$\alpha$ models in domains
with boundaries from this point of view.
\\

In this paper, we investigate the initial boundary value problem
for the LNS-$\alpha$ model (1.1), (1.2) in the following
equivalent form
\begin{eqnarray}\label{1.1}
&& \partial_t v - \Delta v + \nabla\times v \times T_\alpha v +
\nabla p=0 \ {\rm in} \ \Omega\\\label{1.1.1}
&& \nabla\cdot v = 0 \ {\rm in} \ \Omega
\end{eqnarray}
in a smooth bounded domain with the property that both $\Omega$
and $\partial\Omega$ have only finite many simply connected components, where $\nabla\cdot$ and $\nabla\times$
denote the div and curl operator, respectively.
\\

Once the filter mapping $T_\alpha$ is given, equations (1.5) and (1.6)
become a Navier-Stokes type system for $v$, and for which,
some boundary conditions are needed.
Here we consider the following vorticity-slip boundary condition (VSB):
\begin{eqnarray}\label{1.2}
  v \cdot n = 0,\  n\times \nabla\times v = \beta v \ {\rm on} \
  \partial\Omega
\end{eqnarray}
Since there is a boundary, the filter $u = T_\alpha v$ can not be determined
by solving (1.3) and (1.4).
Some boundary conditions are also needed. We propose that the filter $u = T_\alpha v$ be determined by
solving the following
Stokes boundary value problem
\begin{eqnarray}
&& u - \alpha \Delta u +
\nabla \tilde{p} = v \ {\rm in} \ \Omega\\
&& \nabla\cdot u = 0 \ {\rm in} \ \Omega
\end{eqnarray}
with the VSB:
\begin{eqnarray}
  u \cdot n = 0,\  n\times \nabla\times u = \beta u  \ {\rm on} \ \partial\Omega
\end{eqnarray}
We also consider the associated boundary value problem for
the Navier-Stokes equations
\begin{eqnarray}\label{1.3}
&& \partial_t v^0  -\nu\Delta v^0  + (\nabla\times v^0 )\times v^0+
\nabla p = 0\\
&& \nabla\cdot v^0 = 0\label{1.4}
\end{eqnarray}
with the corresponding boundary conditions (1.7) as a singular limit
problem by passing to the vanishing $\alpha$ limit in (1.5)-(1.10).
\\

The equivalence between (1.1) and (1.5) follows from the vector
formula
\begin{eqnarray}
\nabla (u\cdot v) = u\cdot \nabla v + \nabla u^T \cdot v -
\nabla\times v \times u
\end{eqnarray}
for any divergence free vectors u and v.
\\

There have been extensive studies of the Navier-Stokes systems on bounded
domains with various boundary conditions, such as the well known
no-slip condition and various slip boundary conditions. In
particular, substantial understanding has been achieved for the
well-posedness of initial boundary value problems for the
Navier-Stokes system with these boundary conditions and problems
of vanishing viscosity limit and boundary layers, see
\cite{BF,BF1,BdV,BN10,BN04,CQ,CL,C88,GA,IF,KA,Li,TT,YG} and the
references therein. Note that the no-slip boundary condition
corresponds to our VSB with $\beta = \infty$. Yet one of the main
motivations for the proposed VSB is its relation to the well known
Navier-slip boundary condition (see
\cite{Ap,BF,BdV,B67,ITT,IF,So,YG} and the references therein).
Indeed, the Navier-slip boundary condition (NSB) says that the
fluid at the boundary is allowed to slip and the slip velocity is
proportional to the shear stress (see \cite{Na}), i.e.,
\begin{equation}
   v\cdot n = 0, \ 2((S(v)n) )_\tau = -\gamma v_\tau \ {\rm on} \ \partial\Omega
    \end{equation}
where $ 2S(v) = (\nabla v + (\nabla v)^T)$ is the stress tensor. Note that
\begin{equation}
   (2(S(v)n) - (\nabla\times v)\times n)_\tau = GD(v)_\tau \ {\rm on} \ \partial\Omega
    \end{equation}
where $GD(v) = -2S(n)v$ is the lower order term due to the geometry of the boundary, see lemma 3.10. below.
In the special case that the boundary $\partial\Omega$ is flat, one has
$GD(v) = 0$. Thus the VSB (1.7) coincides with NSB (1.14). It should be mentioned
that as far as we know, all the previous physical and numerical studies concerning
the NSB deal with only the case of flat boundaries \cite{Ap,B67,Ja,Jo,PT,TT}. Another
main motivation for the proposed VSB (1.7) and (1.10) is that the vorticity formulations
of the fluid equations have played important roles in analyzing fluid motions, and
suitable boundary conditions on the vorticity should be important for such formulations,
see \cite{Bei,BN10,BN04,CQ,C01,MB} and the references therein. For example, the equivalent
vorticity form of the NSB conditions are crucial in the studies of the corresponding
boundary value problems in \cite{BR,CL}, and the VSB (with $\beta = 0$) was found very
useful to understand the vanishing viscosity limit problem of the Navier-Stokes equations
in \cite{BF, Li, XX, XX2}. It is hoped that the VSB conditions proposed here can share light
on understanding the fluid motions in bounded domains.
\\

The rest of the paper is organized as follows: First, as a
preparation, we present in the next section a $L^2$ version of the
general Hodge decomposition theory that was stated in \cite{CDe} for
smooth vector spaces, which will be used to study the Stokes problems
associated with various slip boundary conditions. Then we give
general and systematic results on well-posedness and spectrum
properties of the Stokes operators associated with various VSB and NSB
conditions in section 3. Our results apply to domains with general
topology. It should be mentioned that all the previous analysis
deals with only the NSB conditions in some special cases. Based on
the properties of the Stokes operators, in section 4, we can formulate
the initial boundary value problem of the LNS-$\alpha$ model, (1.5)-(1.10),
together with the limit problem of the NS equations, (1.11),(1.12),(1.7),
to be a series of abstract equations in a Hilbert space for the
parameter $\alpha\in [0,\infty)$. In section 5, we study the
well-posedness of the weak solutions for the LNS-$\alpha$
equations with the VSB conditions for each $\alpha> 0$, by the
Galerkin method. The local well-posedness, theorem 5.1., is
obtained by direct estimates on the velocity $v$, while the global
theory, theorem 5.2., is proved by combining energy estimates on both
the velocity field $v$ and the filter $u$. Note that our approach is
somewhat different from \cite{FHT,MS} in emphasizing the velocity
$v$ but not the filter $u$. In section 6, we investigate the
vanishing $\alpha$ limit of solutions of the initial boundary value
problem of the LNS-$\alpha$ equations with VSB condition to the
corresponding solutions of the NS equations.  The global in time
convergence of weak solutions is obtained in theorem 6.1. similar
to periodic case in \cite{FHT}, while local in time convergence of
strong solutions is given in theorem 6.2.. The existence of the
global weak solutions and local unique strong solution for the NS
equations with corresponding VSB condition are then followed.
Furthermore, some estimates on convergence rates are given in
theorem 6.3.. Finally, we present some generalizations in section
7. In particular, a parallel theory holds for the NSB condition.

\section{Preparations and Hodge decompositions}

The Hodge decomposition theory plays an important role in the analysis
of vector spaces in a 3D bounded smooth domain, our analysis
on the boundary conditions will be based on this theory. To be self content,
we give a simple $L^2$ version below. For more
details, we refer \cite{CDe, Sc} and the references therein.\\

Let $\Omega\subset R^3$ be a bounded smooth domain, $H^s(\Omega)$ denote the standard Hilbert
space with $H^0(\Omega)= L^2(\Omega)$. Then the
following estimate is well known.
\begin{eqnarray}\label{2.1}
  \|v\|_s\leq c (\|\nabla\times v\|_{s-1} + \|\nabla\cdot v\|_{s-1}
  + |n\cdot v|_{s-\frac{1}{2}} + \|v\|_{s-1})
\end{eqnarray}
for all $v\in H^s(\Omega),\ s\geq 1$ (see \cite{BB, FT}).\\

Let $u\in L^2(\Omega)$. Set
\[ u = v + \nabla \varphi_g \]
Note that $\nabla\cdot u \in\ H^{-1}(\Omega)$. Let $\varphi_g$
solve
\begin{eqnarray}\label{2.2}
&& \Delta \varphi_g= \nabla\cdot u\ {\rm in} \ \Omega\\\label{2.3}
&&  \varphi_g = 0\ {\rm on} \
\partial\Omega
\end{eqnarray}
It follows that
\begin{eqnarray}
  \nabla\cdot v =0\ {\rm in} \ \Omega
\end{eqnarray}
Set
\[DF = \{u\in L^2(\Omega);  \nabla \cdot u = 0\}\]
\[GG = \{u\in L^2(\Omega); u = \nabla \varphi, \varphi\in H_0^1(\Omega)\}\]
Note that
\begin{eqnarray}
  (u, \nabla \varphi) = 0,\ \forall u\in DF,  \varphi\in H_0^1(\Omega)
  \end{eqnarray}
One has

\begin{lem} The following decomposition holds:
\begin{eqnarray}
  L^2(\Omega) = DF \oplus GG
\end{eqnarray}
\end{lem}
Let $u\in DF$. Then $u\cdot n$ is well-defined on $\partial\Omega$
(see \cite{GA}) and
\begin{eqnarray}
 \int_{\partial\Omega}u\cdot n = \int_\Omega \nabla \cdot u  = 0
\end{eqnarray}
Let $\varphi$ solve
\begin{eqnarray} \label{2.4}
&& \Delta \varphi= 0\ {\rm in} \ \Omega\\\label{2.5}
&&  \partial_n\varphi = u\cdot n \ {\rm on} \
\partial\Omega
\end{eqnarray}
Set
\[ v = u - \nabla \varphi \]
and
\[H = \{u\in L^2(\Omega); \nabla \cdot u = 0,  {\rm in} \ \Omega; \ u\cdot n = 0; \ {\rm on}\ \partial\Omega\}\]
\[DFG = \{u\in L^2(\Omega); u = \nabla \varphi, \nabla \cdot u = 0, \int_{\partial\Omega}u\cdot n =0\}\]
Hence $v\in H$.
Note that
\begin{eqnarray}
  (u, v) = 0,\ \forall u\in H, v\in DFG
  \end{eqnarray}
It then follows that

\begin{lem} $DF$ has the following decomposition:
\begin{eqnarray}
   DF = H \oplus DFG
\end{eqnarray}
\end{lem}
Note that $u =\nabla \varphi\in DFG$ may not belong to the range of curl, and
the range of curl $\nabla\times H^1(\Omega)$ is closed in $L^2(\Omega)$. $DFG$ can be further
decomposed to
\begin{eqnarray}
   DFG =  CG \oplus HG
\end{eqnarray}
where
\[  CG =  DFG\cap (\nabla\times H^1(\Omega)), \ HG = DFG\cap (\nabla\times H^1(\Omega))^\perp \]
\\
Let $u= \nabla \varphi\in HG$. Since
\begin{eqnarray}
 0 = ((\nabla\times v),\nabla \varphi)
  = \int_{\partial\Omega} (n\times v)\cdot \nabla \varphi
\end{eqnarray}
for all $v\in H^1(\Omega)$, thus $\partial_\tau\varphi  = 0,\ {\rm on} \
\partial\Omega$ with $\tau$ being any tangential direction on $\partial\Omega$
which implies $\varphi$ is a constant on each component $\Gamma_i$ of
$\partial\Omega$. So
\[HG  = \{\nabla \varphi; \Delta\varphi = 0, \varphi = c_i\ on\ \ \Gamma_i\}\]
consists only smooth vectors, and is finite dimensional, which is
called the harmonic gradient space.

\begin{rem} $CG$ can also be expressed as
\[CG = \{u\in L^2(\Omega);\ u = \nabla \varphi, \nabla\cdot u = 0,\ \int_{\Gamma_i}u\cdot n = 0\}\]
\\
Since $CG\subset \nabla\times H^1(\Omega)$,
we will call it curl type gradient space.
\end{rem}
Note that $H\cap Ker(\nabla\times)$ is compact in $L^2(\Omega)$ due to (2.1).
Set
\[ HH = H\cap Ker(\nabla\times)\]
Then
\[ HH =
\{u\in L^2(\Omega); \nabla \cdot u = 0,\ \nabla\times u = 0\ {\rm
in} \ \Omega, \ u\cdot n = 0; \ {\rm on}\ \partial\Omega\}\] This
is called the harmonic knots space, which consists only smooth
functions and is finite dimensional (see \cite{CDe}). Now $H$ can
be decomposed to
\begin{eqnarray}
   H = FH \oplus HH
\end{eqnarray}
where
\[ FH =  H\cap (Ker(\nabla\times))^\bot\]
\\
In conclusion, we have

\begin{lem} The following decomposition holds:
\begin{eqnarray}
  L^2(\Omega) =  FH \oplus HH \oplus CG \oplus HG \oplus GG
\end{eqnarray}
\end{lem}
Then for any $u\in  L^2(\Omega)$, it is uniquely written to
\begin{eqnarray}
  u = P_{HH} u + P_{FH} u  + P_{CG}u  + P_{HG}u + P_{GG} u
\end{eqnarray}
where $P_{X}$ denotes the projection on the corresponding subspace.\\
\\
It should be noticed that the space $FH$ has the following expressions (see [9, 20, 36, 47]).

\begin{lem} The space $FH$ can be expressed as
\begin{eqnarray}
 FH = \{u\in L^2(\Omega); \nabla\cdot u = 0, u\cdot n = 0\ on\ \partial\Omega, F(u)=0 \}
 \end{eqnarray}
 \begin{eqnarray}
 FH = \{u; u = \nabla\times v,\ v\in H^1(\Omega),\ \nabla\cdot v= 0, \ n\times v = 0\ on\ \partial\Omega\}
 \end{eqnarray}
where $F(u)=0$ means
\[   \int_\Sigma u\cdot n = 0\]
for any smooth cross section $\Sigma$ of $\Omega$.
\end{lem}
It follow from (2.1) and the fact that $HH\subset \nabla\times (FH\cap H^1(\Omega))$
(see \cite{CDe}) respectively that

\begin{prop}
\begin{eqnarray}
  L^2(\Omega) = \nabla\times (FH\cap H^1(\Omega))\oplus HG \oplus GG
\end{eqnarray}
\end{prop}
Similarly, in general, it holds that

\begin{prop}
\begin{eqnarray}
H^s(\Omega) = \nabla\times (FH\cap H^{s+1}(\Omega))\oplus (HG\cap H^s(\Omega)) \oplus (GG \cap H^s(\Omega))
\end{eqnarray}
for $s\geq 0$.
\end{prop}
It follows from (\ref{2.1}),(\ref{2.2}),(\ref{2.3}),(\ref{2.4}),(\ref{2.5}) and the fact
that $HH$, $HG$ are finite dimensional that

\begin{prop} $C^\infty(\Omega)\cap X$ is dense in
$H^s(\Omega)\cap X,\ s\geq 0$ for
\[X = FH, HH, CG, HG, GG \]
\end{prop}

\section{The Stokes operators}

In this section, we apply the Hodge decomposition theory to the
Stokes problems with both the VSB and NSB conditions. We first
consider a special Stokes problem with the VSB (3.1)-(3.3) and
prove theorem 3.1.. Next, since the topology of the domain is
assumed to be general, to avoid the uniqueness of the solutions
for the general Stokes problems, we consider the perturbed Stokes
problem associated with VSB (\ref{3.1})-(\ref{3.3}). Based on
theorem 3.1., by using the Hodge decomposition theory, we prove
the associated Stokes operator is a self-adjoint extension of the
associated positive definite bilinear form (see theorem 3.5.). The
proof of theorem 3.5. is constructive, and the techniques can also be
used to prove the well-posedness of the non-homogeneous problem
(\ref{3.6})-(\ref{3.7})(see theorem 3.7.). More generally, we can
prove well-posedness of the boundary value problem
(\ref{3.8})-(\ref{3.9}) (see theorem 3.9.) by construction a
contraction map. Finally, we identify the relationship between the NSB
and VSB, and establish a similar theory for the Stokes
problem associated with NSB (\ref{3.1n})-(\ref{3.3n}).

\subsection{A special Stokes problem}

Let us start by considering the following special Stokes problem with homogenous VSB condition
\begin{eqnarray}
&& - \Delta u  = f \ {\rm in} \ \Omega\\
&& \nabla\cdot u = 0 \ {\rm in} \ \Omega\\
&& u \cdot n = 0, \  n\times \nabla\times u = 0 \ {\rm on} \ \partial\Omega
\end{eqnarray}
with $f\in FH$. Set
\[ W  = \{u\in H^2(\Omega); \ n\times (\nabla\times u) = 0 \ {\rm on},
\partial\Omega\}\]
Then we have

\begin{thm} The Stokes operator $A_F = -\Delta$ with the domain
$D(A_F) = W \cap FH$ is self-adjoint in the Hilbert space $FH$.
\end{thm}
\textbf{Proof:}
It is clear that $A_F = -\Delta$ with the domain $W \cap FH$
is symmetric. Since $C_0^\infty(\Omega)\cap H$ is dense in $H$,
it follows that $A_F$ is densely defined due to the orthogonality of $FH$ and $HH$ and
the compactness of $HH$. Let $u\in W$.
Since $n\times(\nabla\times u) = 0$ on $\partial\Omega$, then
$-\Delta u =\nabla\times (\nabla\times u) \in FH$ by lemma 2.5., thus
$A_F$ maps $W \cap FH$ to $FH$. Now, for any $f\in FH$,
it follows from lemma 2.5. that
there is a $\Phi\in H^1(\Omega)$ satisfying
\begin{eqnarray}
&&  \nabla\times \Phi = f \ {\rm in} \ \Omega\\
&&  \nabla\cdot \Phi = 0 \ {\rm in} \ \Omega\\
&&  \Phi\times n = 0 \ {\rm on} \ \partial\Omega
\end{eqnarray}
Due to proposition 2.7. and lemma 2.1., there is a $v\in FH\cap H^2(\Omega)$ so that
\begin{eqnarray}\label{3.1.1}
  \Phi = \nabla\times v + P_{HG}\Phi
\end{eqnarray}
Note that $P_{HG}\Phi\times n = 0$ on $\partial\Omega$. It
follows that
\begin{eqnarray}
  n\times(\nabla\times v) = 0\ on\ \partial\Omega
\end{eqnarray}
Then $\nabla\times(P_{HG}\Phi)=0$ and (\ref{3.1.1}) imply that
\begin{eqnarray}
 -\Delta v = f\ {\rm in} \ \Omega
\end{eqnarray}
Thus $A_F: W \cap FH\rightarrow FH$ is surjective.
If $f=0$, then integration by part shows
\begin{eqnarray}
  \|\nabla\times v\| = 0
\end{eqnarray}
It follows that $u = 0$ due to the orthogonality of $FH$ and $HH$ and then
$A_F: W \cap FH\rightarrow FH$ is one to one.\\
\\
Noting that $W$ and $FH$ are closed in $H^2(\Omega)$ and
$L^2(\Omega)$, and
\begin{eqnarray}
 \|\Delta v\| \leq \|v\|_2
\end{eqnarray}
we obtain from the Banach inverse operator theorem that
\begin{eqnarray}
 \|v\|_2 \leq c\|\Delta v\|
\end{eqnarray}
The theorem was proved.
\\

Equivalently, we have shown the problem (3.1)-(3.3)
has a unique solution $u\in H^2(\Omega)$ for any $f\in FH$. \\

It follows from the proof of theorem 3.1. that
\[ \nabla\times: H_n^1(\Omega)  \mapsto HF\]
is also one to one and onto, where
\[ H_n^1(\Omega)= \{u\in  H^1(\Omega); \nabla\cdot u =0; \ \ n\times u =0\ on\ \partial \Omega; (u,\varphi)=0, \forall \varphi\in HG\}\]
It follows from the trace theorem and the continuity of the divergence
operator that $H_n^1(\Omega)$ is closed in
$H^1(\Omega)$. Then
\begin{eqnarray}
 \|u\|_1 \leq c\|\nabla\times u\|
\end{eqnarray}
follows from
\begin{eqnarray}
 \|\nabla\times u\|\leq \|u\|_1
\end{eqnarray}
for any $u\in H_n^1(\Omega)$. This yields immediately that

\begin{lem} Let  $u \in H_n^1(\Omega)$. Then the following
Poincar\'{e} type inequality holds
\begin{eqnarray}
 \|u\|\leq c\|\nabla\times u\|
\end{eqnarray}
\end{lem}
Let $v\in FH\cap H^1(\Omega)$. Then there is $u\in H_n^1(\Omega)$ such that $\nabla\times u = v$
and
\begin{eqnarray}
 (v,v) = (\nabla\times u, v) = (u,\nabla\times v) \leq \|u\|\|\nabla\times v\|
\end{eqnarray}
This, together with (3.15), shows
\begin{eqnarray}
 (v,v) \leq c\|\nabla\times v\|\|\nabla\times u\|
\end{eqnarray}
Thus, one gets

\begin{lem} Let  $u\in FH\cap H^1(\Omega)$. Then the
following Poincar\'{e} type inequality
\begin{eqnarray}
 \|u\|\leq c\|\nabla\times u\|
\end{eqnarray}
holds.
\end{lem}
As a consequence, we can obtain

\begin{cor} The operator $A_F$ in theorem 3.1. is the
self adjoint extension of the following bilinear form
\begin{eqnarray}
 a(u,\phi) = (\nabla\times u, \nabla\times \phi), \ \emph{D}(a)= V_F=FH\cap H^1(\Omega)
\end{eqnarray}
in $FH$.
\end{cor}
\textbf{Proof:} From proposition 2.8., $a(u,\phi)$ with $\emph{D}(a)= FH\cap H^1(\Omega)$
is densely defined. Due to (2.1) and lemma 3.3., $a(u,\phi)$ is closed and positive.
It follows that there is a self-adjoint operator $A$ with domain $\emph{D}(A)\subset\emph{D}(a)$ such that
\begin{eqnarray}
 a(u,\phi) = (Au, \phi), \forall\phi\in FH\cap H^1(\Omega)
\end{eqnarray}
for any $u\in \emph{D}(A)$. It is clear that $\emph{D}(A_F) =
W\cap FH\subset \emph{D}(A)$ and $Au = -\Delta u$ for any $u\in
W\cap FH$. Let $u\in \emph{D}(A)$ and $f = Au$. It then follows
that $f\in FH$. It follows from theorem 3.1. that there is a $v\in
\emph{D}(A_F)$ such that (3.1)-(3.3) are valid (with $u$ replaced by $v$) and
\begin{eqnarray}
 a(v,\phi) = (f,\phi)
\end{eqnarray}
for all $\phi\in V_F$. On the other hand
\begin{eqnarray}
 a(u,\phi) = (Au, \phi) = (f,\phi)
\end{eqnarray}
for all $\phi\in V_F$, hence
\begin{eqnarray}
 a(u-v,\phi) = (\nabla\times (u-v), \nabla\times \phi) =0
\end{eqnarray}
for all $\phi\in V_F$. Taking $\phi= u-v$ shows that
$\nabla\times (u-v) = 0$. Thus $u=v$ due to (2.14).
Thus $\emph{D}(A) = \emph{D}(A_F)$ and $A= A_F$.
\\

Denote by $V_F'$ the dual space of $V_F$ respect to the $L^2$ inner product. Then the
notation of weak solutions can be extended for $f\in V_F'$: $u$ is called a weak solution to (3.1)-(3.3) for $f\in V_F'$ if
\begin{eqnarray}
 a(u,\phi) = (f,\phi),\ \forall\ \phi\in V_F
\end{eqnarray}

\subsection{The Stokes problem with VSB condition}
Next, we consider the Stokes problem with general VSB condition.
Since the domain is allowed to have general topologe, the kernel
of $-\Delta$ may be not empty. To avoid it, we consider the
following boundary value problem instead:
\begin{eqnarray}\label{3.1}
&&  (I- \Delta) u + \nabla p = f \ {\rm in} \ \Omega\\ \label{3.2}
&& \nabla\cdot u = 0 \ {\rm in} \ \Omega\\ \label{3.3} && u\cdot
n = 0,\  n\times(\nabla\times u) = \beta u \ {\rm on} \
\partial\Omega
\end{eqnarray}
where $\beta$ is a nonnegative smooth function.
\\
Define
\[ V = H^1(\Omega)\cap H \]
\[ W_\beta = \{u\in H^2(\Omega); \ n\times (\nabla\times u) = \beta u \ {\rm on},
\partial\Omega\}\]
Define a bilinear form as
\[ \tilde{a}_\beta(u,\phi) =  (u,\phi)+ a_\beta(u,\phi) \]
where
\begin{eqnarray}
a_\beta(u,\phi) = \int_{\partial\Omega}\beta u\cdot \phi + \int_\Omega
(\nabla\times u)\cdot (\nabla\times \phi)
\end{eqnarray}
with the domain $\emph{D}(\tilde{a}_\beta)= V$. $u\in V$ is said to be a weak solution to
the boundary value problem (\ref{3.1})-(\ref{3.3}) on $H$ for $f\in V'$ if
\begin{eqnarray}
  \tilde{a}_\beta(u,\phi) = (f, \phi), \  \forall \ \phi\ {\rm in}\ V
\end{eqnarray}
where $V'$ is the dual space of $V$. Based on theorem 3.1., we can prove

\begin{thm}The self-adjoint extension of the bilinear form
$\tilde{a}_\beta(u,\phi)$ with the domain $\emph{D}(\tilde{a}_\beta) = V $ is the
Stokes operator $A_\beta =  I+ P(-\Delta)$ with $\emph{D}(A_\beta) =
W_\beta\cap H $, and $A_\beta$ is an isomorphism between $\emph{D}(A_\beta)$ and
$H$ with compact inverse on $H$. Consequently, the
eigenvalues of the Stokes operator $A_\beta$ can be listed as
\[ 1 \leq 1+\lambda_{1}\leq 1+ \lambda_{2}\cdots \rightarrow \infty \]
with the corresponding eigenvectors $\{e_{j}\}\subset W_\beta$, i.e.,
\begin{eqnarray}\label{3.4}
   A_\beta e_j = (1+ \lambda_j)e_j
\end{eqnarray}
which form a complete orthogonal basis in $H$. Furthermore, it
holds that
\begin{eqnarray}\label{3.5}
   (1+\lambda_1)\|u\|^2\leq \tilde{a}_\beta(u,u) \leq \frac{1}{1+ \lambda_1} \|A_\beta u\|^2,\
   \forall\ u\in \emph{D}(A_\beta)
\end{eqnarray}
\end{thm}
\textbf{Proof:} It is clear that $\tilde{a}_\beta(u,\phi)$ with the domain
$\emph{D}(\tilde{a}_\beta) = V $ is a positive densely defined
closed bilinear form. Let $A_\beta$ be the self-adjoint extension of
$\tilde{a}_\beta(u,\phi)$. It follows that $ W_\beta\cap H \subset \emph{D}(A_\beta)$
and $A_\beta u = u + P(-\Delta u)$, for any $u\in W_\beta\cap H$ by integrating by part.
It remains to show that $\emph{D}(A_\beta)\subset W_\beta\cap H$.
Let $u\in \emph{D}(A_\beta)$ and $f = A_\beta u$. Since $\emph{D}(A_\beta)\subset
\emph{D}(\tilde{a}_\beta) = V$, it follows from (3.29) that
\begin{eqnarray}\label{3.11}
  \|u\|_{1} \leq  c \|f\|
\end{eqnarray}
Let $n(x)$ and $\beta(x)$ be internal smooth extensions of the
normal vector $\beta$ respectively. Then $\beta(x)u\times n(x)\in
H^1(\Omega)$. Proposition 2.7. yields
\begin{eqnarray}
 \beta(x)u\times n(x) = \nabla\times v + \nabla h + \nabla g
\end{eqnarray}
with $\nabla h = P_{HG}(\beta(x)u\times n(x))$,
$\nabla g =P_{GG}(\beta(x)u\times n(x))$ and $v\in FH\cap H^2(\Omega)$. It follows that
\begin{eqnarray}
\|\nabla g\|_1\leq c \|u\|_1
\end{eqnarray}
since $g$ satisfies
\begin{eqnarray}
&& \Delta g = \nabla\cdot(\beta(x)u\times n(x))\ {\rm in} \ \Omega\\
&&   g = 0\ {\rm on} \ \partial\Omega
\end{eqnarray}
Since $HG$ is finite dimensional, so
\begin{eqnarray}
\|P_{HG}(\beta(x)u\times n(x))\|_1 \leq c\|P_{HG}(\beta(x)u\times n(x))\|\leq c\|u\|
\end{eqnarray}
It then follows from (2.1) and lemma 3.3 that
\begin{eqnarray}\label{3.12}
&& \|v\|_2 \leq c \|\nabla\times v\|_1 \leq c\|u\|_1\leq  c \|f\|
\end{eqnarray}
Integrating by part and noting that $n\times \nabla h = 0$, $n\times \nabla g = 0$
on the boundary, we have
\begin{eqnarray}
 \int_\Omega (\nabla\times v)\cdot (\nabla\times \phi) +
 \int_{\partial\Omega}\beta n\times(u\times n)\cdot \phi = (-\Delta v,\phi)
\end{eqnarray}
for all $\phi\in H^1(\Omega)$.
It follows from $n\times(u\times n)= u$ on the boundary and the definition of the weak solution that
\begin{eqnarray}
 \int_\Omega (\nabla\times (u-v))\cdot (\nabla\times \phi)=
 (P_{FH}(f-u+\Delta v), \phi)
\end{eqnarray}
for all $\phi\in H^1(\Omega)\cap FH$.
Note that $\nabla\times u = \nabla\times P_{FH}(u)$ and $P_{FH}(u)\in H^1(\Omega)\cap FH$.
It follows that
\begin{eqnarray}
a(P_{FH}(u)-v,\phi) = (P_{FH}(f-u+\Delta v), \phi),\ \forall \phi\ \in \ H^1(\Omega)\cap FH
\end{eqnarray}
It follows from theorem 3.1. that $P_{FH}(u)-v \in W$ and
\begin{eqnarray}\label{3.13}
 \|P_{FH}(u)-v\|_2\leq c(\|f\|+ \|\Delta v\|+ \|u\|)
\end{eqnarray}
Since $HH$ is finite dimensional, it holds that
\begin{eqnarray}\label{3.14}
 \|P_{HH}(u)\|_2\leq c \|u\|
\end{eqnarray}
One gets from (\ref{3.11}),(\ref{3.12}),(\ref{3.13}) and (\ref{3.14}) that
\begin{eqnarray}
 \|u\|_2\leq c\|f\|
\end{eqnarray}
Since $P_{FH}(u)-v \in W$, it holds that
\begin{eqnarray}
 n\times\nabla\times u = n\times\nabla\times P_{FH}(u)= n\times\nabla\times v = \beta u \ {\rm on} \ \partial\Omega
\end{eqnarray}
Thus we have shown $u\in W_\beta\cap H$. Integrating by part in (3.29) yields
\begin{eqnarray}
 (u-\Delta u-f, \phi) = 0
\end{eqnarray}
for all $\phi\in V$, which implies
\begin{eqnarray}
&& u - \Delta u + \nabla p = f \ {\rm in} \ \Omega
\end{eqnarray}
with $p$ given by
\begin{eqnarray}
&& -\Delta p = 0\ {\rm in} \ \Omega\\
&& (\nabla p)\cdot n =  \Delta u\cdot n \ {\rm on} \
\partial\Omega
\end{eqnarray}
It is noted that $\|A_\beta u\|$ is an equivalent norm of
$H^2(\Omega)$ on $ W_\beta\cap H$ due to (3.44) and
\begin{eqnarray}
 \|A_\beta u\|\leq\|u\|+ \|\Delta u\|\leq c \|u\|_2
\end{eqnarray}
for all $u\in W_\beta\cap H$. The theorem was proved.
\\

Let $V'$ be the dual space of $V$ with respect to the $L^2$ inner
product. $u\in V$ is called a weak solution to (3.1)-(3.3) for $f\in
V'$ if
\begin{eqnarray}
 \tilde{a}_\beta(u,\phi) = (f,\phi),\ \forall\ \phi\in V
\end{eqnarray}
By using a standard density argument, one can show

\begin{cor} For any $f\in V'$, the boundary value
problem (\ref{3.1})-(\ref{3.3}) has a unique weak solution $u\in
V$
\end{cor}
Now, let $b\in H^{\frac{1}{2}}(\partial\Omega)$ and $b\cdot n =
0$ on $\partial\Omega$. From the extension theorem, it has an extension denoted by
$b(x)\in H^1(\Omega)$. Similar to the proof of theorem 3.5., one can show that there
exists a $\Phi\in H^2(\Omega)\cap FH$ such that
\[ n\times(\nabla\times\Phi) = b \ {\rm on} \ \partial\Omega \]
It follows that $\Phi$ solves the following problem:
\begin{eqnarray}
&& u - \Delta u + \nabla p = f \ {\rm in} \ \Omega\\
&& \nabla\cdot u = 0 \ {\rm in} \ \Omega\\
&& u\cdot n = 0, \ n\times(\nabla\times u) = b\ {\rm on} \
\partial\Omega
\end{eqnarray}
with $f = u + P(-\Delta \Phi)$ and $\nabla p = \Delta \Phi -
P(\Delta \Phi)$. This fact and theorem 3.5. for $\beta = 0$ yield

\begin{thm} Let $b\in H^{\frac{1}{2}}(\partial\Omega)$,
$b\cdot n = 0$ and $\lambda
>0$. Then the following problem
\begin{eqnarray}\label {3.6}
&& \lambda u - \Delta u + \nabla p = f \ {\rm in} \ \Omega\\
&& \nabla\cdot u = 0 \ {\rm in} \ \Omega\\\label {3.7}
&& u\cdot
n = 0, \ n\times(\nabla\times u) = b \ {\rm on} \
\partial\Omega
\end{eqnarray}
has a unique solution $u\in H^2(\Omega)$ for any $f\in H$.
\end{thm}
The boundary value problem (\ref{3.6})-(\ref{3.7}) also have a
weak formulation
\begin{eqnarray}\label{3.10b}
  \lambda(u,\phi)+ \int_\Omega (\nabla\times u)\cdot(\nabla\times \phi) + \int_{\partial\Omega}b\cdot \phi
   = (f, \phi), \  \forall \ \phi\ {\rm in}\ V
\end{eqnarray}
Similar to corollary 3.4., one has

\begin{cor} Let $b\in
H^{-\frac{1}{2}}(\partial\Omega)$, $b\cdot n = 0$. Then for any
$f\in V'$, the boundary value problem (\ref{3.6})-(\ref{3.7}) has
a unique weak solution $u\in V$
in the sense of (\ref{3.10b}).
\end{cor}
We omit the details of the proof here, and refer to \cite{GR} for the
definition of the weak tangential trace $H^{-\frac{1}{2}}(\partial\Omega)$.
\\

For any given smooth and nonnegative function $\beta$, we define
the map
\[ T : H^{\frac{1}{2}}(\Omega)\cap H \mapsto V\subset  H^{\frac{1}{2}}(\Omega)\cap H \]
by $ u = Tv $ determined by (\ref{3.10b}) with $b$ replaced by
$\beta v + b$ and $ f= 0$. Let $v_i \in H^{\frac{1}{2}}(\Omega)\cap H$ and $u_i=Tv_i$, $i=1,2$. It then follows from (3.58) that
\begin{eqnarray}
 \lambda\|u_1- u_2\|^2 + \|\nabla\times(u_1- u_2)\|^2
 + \int_{\partial\Omega}\beta (u_1- u_2)\cdot(v_1- v_2) = 0
\end{eqnarray}
Note that
\begin{eqnarray}
 |\int_{\partial\Omega}\beta (u_1- u_2)\cdot(v_1- v_2) |\leq
 c\|u_1-u_2\|_{H^{\frac{1}{2}}(\Omega)}\|v_1-v_2\|_{H^{\frac{1}{2}}(\Omega)}
\end{eqnarray}
and
\begin{eqnarray}
 \|\varphi\|_{H^{\frac{1}{2}}(\Omega)}^2\leq
 c\|\varphi\|\|\varphi\|_{H^1(\Omega)}\leq
 c\|\varphi\|\|\nabla\times\varphi\|,\ \forall\ \varphi\in V
\end{eqnarray}
It follows that
\begin{eqnarray}
 \|u_1-u_2\|_{H^{\frac{1}{2}}(\Omega)}^2\leq
 c\lambda^{-\frac{1}{2}}\|v_1-v_2\|_{H^{\frac{1}{2}}(\Omega)}^2
\end{eqnarray}
for $\lambda \geq 1$. Take $\lambda$ large enough such that $T$
becomes a contraction map on $H^{\frac{1}{2}}(\Omega)$. It follows
that
\[  Tv = v \]
has a unique solution $\Psi$ on $H^{\frac{1}{2}}(\Omega)$ and
$\Psi = T\Psi \in H^1(\Omega)$. \\

For any $\tilde{f}\in V'$, let $v$ be the weak solution of
(\ref{3.1})-(\ref{3.3}) with $f = \tilde{f} -(1-\lambda)\Psi$. It is
clear that $u = v + \Psi\in V$ is a weak solution of the following
problem:
\begin{eqnarray}\label {3.8}
&& u - \Delta u + \nabla p = \tilde{f} \ {\rm in} \ \Omega\\
&& \nabla\cdot u = 0 \ {\rm in} \ \Omega\\\label {3.9}
&& u\cdot
n = 0, \ n\times(\nabla\times u) = \beta u + b \ {\rm on} \
\partial\Omega
\end{eqnarray}
in the sense that
\begin{eqnarray}\label{3.10}
  (u,\phi)+ \int_\Omega (\nabla\times u)\cdot(\nabla\times \phi) + \int_{\partial\Omega}(\beta u+b)\cdot \phi
   = (\tilde{f}, \phi), \  \forall \ \phi\ {\rm in}\ V
\end{eqnarray}
The uniqueness can be proved in the same way as for theorem 3.1.. We conclude

\begin{thm} Let $b\in H^{-\frac{1}{2}}(\partial\Omega)$,
$b\cdot n = 0$. Then for any $\tilde{f}\in V'$, the boundary value
problem (\ref{3.8})-(\ref{3.9}) has a unique solution $u\in V$ in
the sense of (\ref{3.10}). Moreover, if $\tilde{f}\in H$ and $b\in
H^{\frac{1}{2}}(\partial\Omega)$,
then $u\in H^2(\Omega)$.
\end{thm}

\subsection{The Stokes problem with the NSB condition}

We can establish a similar theory for the Stokes problem with the NSB
just as with VSB. For completeness, we sketch it here. Consider
the following Stokes problem with the NSB condition.
\begin{eqnarray}\label{3.1n}
&&  (I- \Delta)u + \nabla p = f \ {\rm in} \ \Omega\\\label{3.2n}
&&\nabla\cdot u = 0 \ {\rm in} \ \Omega\\\label{3.3n} && u\cdot n
= 0,\ 2(S(u)n)_{\tau} =-\gamma u_{\tau} u \ {\rm on} \
\partial\Omega
\end{eqnarray}
where $\gamma$ is a nonnegative smooth function.\\
Define
\[ \tilde{W}_\gamma = \{u\in H^2(\Omega); \ 2(S(u)n)_{\tau} =-\gamma u_{\tau} \ {\rm on},
\partial\Omega\}\]
\\
and a bilinear form
\[ \tilde{a}_\gamma(u,\phi) = (u,\phi)+  a_\gamma(u,\phi),\ \emph{D}(\tilde{a}_\gamma) = V \]
where
\begin{eqnarray}
a_\gamma(u,\phi) =  \int_{\partial\Omega}\gamma u\cdot \phi +
2\int_\Omega S(u)\cdot S(\phi)
\end{eqnarray}
and $S(u)\cdot S(\phi)$ denotes the trace of the product of the two matrices.\\

$u$ is said to be a weak solution to the boundary value problem
(\ref{3.1n})-(\ref{3.3n}) on $H$ for $f\in V'$ if
\begin{eqnarray}\label{3.26}
  \tilde{a}_\gamma(u,\phi)  = (f, \phi), \  \forall \ \phi\ {\rm in}\ V
\end{eqnarray}
where $V'$ is the dual space of $V$.\\

To compare it with the VSB case,
we first calculate that
\begin{lem} Let $u\in H^2(\Omega)$ and $u\cdot n = 0$ on
the boundary. It holds that
\begin{equation}
   (2(S(u)n) - \omega\times n)_\tau = GD(u)_\tau
    \end{equation}
with $GD(u) = -2S(n)u$.
\end{lem}
\textbf{Proof:} Note that
\begin{equation}
  \partial_n u = \frac{1}{2}\omega\times n + S(u)n
    \end{equation}
and
\begin{equation}
  \partial_\tau u = \frac{1}{2}\omega\times \tau + S(u)\tau
    \end{equation}
It follows that
\begin{eqnarray}
&& 2(S(u)n)\cdot\tau = \partial_\tau u\cdot n + \partial_n u \cdot \tau \\
&& (n\times \omega)\tau = \partial_\tau u\cdot n -\partial_n u
\cdot \tau
    \end{eqnarray}
and
\begin{equation}
   2(S(u)n)\cdot\tau + (n\times \omega)\tau = 2\partial_\tau u\cdot
   n
    \end{equation}
Note that $u\cdot n = 0$ on the boundary. It follows that
\begin{equation}
   \partial_\tau u\cdot n = - u\cdot\partial_\tau n
    \end{equation}
We conclude that
\begin{equation}
  ( 2S(u)n - \omega\times n)\cdot\tau = -2u\cdot\partial_\tau n
    \end{equation}
Note that
\begin{equation}
  \partial_\tau n = \frac{1}{2}(\nabla\times n)\times \tau +
  S(n)\tau
    \end{equation}
thus
\begin{equation}
  ( 2S(u)n - \omega\times n)\cdot\tau = ((\nabla\times n)\times u)\cdot\tau - 2S(n)u\cdot
  \tau
    \end{equation}
Note that
\begin{equation}
  u\times\tau = \lambda n
    \end{equation}
and
\begin{equation}
  (\nabla\times n)\cdot n = 0
    \end{equation}
on the boundary. It follows that
\begin{equation}
  ( 2S(u)n - \omega\times n)\cdot\tau = - 2S(n)u\cdot \tau
    \end{equation}
Set
\[ GD(u) = -2S(n)u\]
The lemma is proved.
\\

It follows from a simple calculation and by using the density method that
\begin{lem} Let $u\in H^1(\Omega)\cap H$. Then
\begin{eqnarray}
 2\int_\Omega S(u)\cdot S(\phi) = \int_\Omega
(\nabla\times u)\cdot (\nabla\times \phi) + \int_{\partial\Omega}
GD(\phi)\cdot u
\end{eqnarray}
\begin{eqnarray}\label{3.21}
\int_{\partial\Omega} GD(\phi)\cdot u = \int_{\partial\Omega}
GD(u)\cdot\phi
\end{eqnarray}
\end{lem}
As a counterpart of theorem 3.5., we can obtain

\begin{thm} The self-adjoint extension of the bilinear
form $\tilde{a}_\gamma(u,\phi)$ with domain
$\emph{D}(\tilde{a}_\gamma) = V $ is the Stokes operator $A_\gamma
=  I+ P(-\Delta)$ with $\emph{D}(A_\gamma) = \tilde{W}_\gamma\cap
H $, and $A_\gamma$ is an isomorphism between $\emph{D}(A_\gamma)$
and $H$ with a compact inverse on $H$. Consequently, the eigenvalues
of the Stokes operator $A_\gamma$ can be listed as
\[ 1 \leq 1+ \lambda_{1}\leq 1+ \lambda_{2}\cdots \rightarrow \infty\]
with the corresponding eigenvectors $\{e_{j}\}\subset
\tilde{W}_\gamma$, i.e.,
\begin{eqnarray}\label{3.41}
   A_\gamma e_j = (1+ \lambda_j)e_j
\end{eqnarray}
which form a complete orthogonal basis in $H$. Furthermore, it
holds that
\begin{eqnarray}
   (1+\lambda_1)\|u\|^2\leq \tilde{a}_\beta(u,u) \leq \frac{1}{1+\lambda_1} \|A_\gamma
   u\|^2,\ \forall u\in \emph{D}(A_\gamma)
\end{eqnarray}
\end{thm}
\textbf{Proof:} It suffices to show that $\emph{D}(A_\gamma)\subset
\tilde{W}_\gamma\cap H$ since the rest is similar to the proof of
theorem 3.5.. Let $u\in \emph{D}(A_\gamma)$ and $f=A_\gamma u$. Since
$\emph{D}(A_\gamma)\subset \emph{D}(\tilde{a}_\beta) =
H^1(\Omega)\cap H$, it follows from (\ref{3.26}) that
\begin{eqnarray}\label{3.22}
  \|u\|_{1}^{2} \leq  c \|f\|^2
\end{eqnarray}
Let $n(x)$ and $\gamma(x)$ be internal smooth extensions of the
normal vector $n$ and $\gamma$. Then $(\gamma(x)u+GD(u))\times
n(x)\in H^1(\Omega)$. Due to proposition 2.7., one has
\begin{eqnarray}
 (\gamma(x)u+GD(u))\times n(x) = \nabla\times v + \nabla h + \nabla g
\end{eqnarray}
with $v\in H^2(\Omega)\cap FH$, $\nabla h =
P_{HG}((\gamma(x)u+GD(u))\times n(x))$ and $\nabla g =
P_{GG}((\gamma(x)u+GD(u))\times n(x))$. Similar to the proof of
theorem 3.1., one can get
\begin{eqnarray}\label{3.23}
&& \|v\|_2 \leq c \|\nabla\times v\|_1 \leq c\|u\|_1
\end{eqnarray}
Note that $n\times (\nabla h)= 0$ and $n\times (\nabla g) =0$.
Thus
\begin{eqnarray}
 \int_\Omega (\nabla\times v)\cdot (\nabla\times \phi)+
 \int_{\partial\Omega}(\gamma u + GD(u))\cdot\phi = (-\Delta v,
 \phi),\ \forall \phi\in V
\end{eqnarray}
Then the definition of the weak solution and lemma 3.11. imply
\begin{eqnarray}
 \int_\Omega (\nabla\times u)\cdot (\nabla\times \phi)+
 \int_{\partial\Omega}\gamma u \cdot\phi + \int_{\partial\Omega} GD(\phi)\cdot u = (f-u, \phi),
 \ \forall \phi\in V
\end{eqnarray}
Combine them and note (\ref{3.21}) to get
\begin{eqnarray}
 \int_\Omega (\nabla\times (u-v))\cdot (\nabla\times \phi)=
 (P_{FH}(f-u+ \Delta v), \phi),\ \forall \phi\in V
\end{eqnarray}
Note that $\nabla\times u = \nabla\times P_{FH}(u)$ and $P_F(u)\in H^1(\Omega)\cap FH$.
It follows that
\begin{eqnarray}
a(P_{FH}(u)-v,\phi) = (P_{FH}(f-u+ \Delta v), \phi),\ \forall
\phi\ \in \ H^1(\Omega)\cap FH
\end{eqnarray}
Since $P_{FH}(f-u+ \Delta v)\in FH$, so $P_{FH}(u)-v\in W$, and
\begin{eqnarray}\label{3.24}
 \|P_{FH}(u)-v\|_2\leq c(\|f\|+ \|u\|_1)
\end{eqnarray}
Since $HH$ is a finite dimensional, so
\begin{eqnarray}\label{3.25}
 \|P_{HH}(u)\|_2\leq c \|u\|
\end{eqnarray}
It follows from (\ref{3.22}),(\ref{3.23}),(\ref{3.24}) and
(\ref{3.25}) that
\begin{eqnarray}
 \|u\|_2\leq c\|f\|
 \end{eqnarray}
Note that
\begin{eqnarray}
 (\nabla\times u)\times n= (\nabla\times P_{FH}(u))\times n= (\nabla\times v)\times n= -\gamma u -
 GD(u)
\end{eqnarray}
It follows that
\begin{eqnarray}
  2(S(u)n)_\tau = ((\nabla\times u)\times n + GD(u))_\tau = -\gamma
  u_\tau
\end{eqnarray}
The theorem was proved.
\\

Similar to the discussion for VSB, we have

\begin{thm} Let $b\in H^{-\frac{1}{2}}(\partial\Omega)$,
$b\cdot n = 0$, $\gamma$ be a nonnegative smooth function on the
boundary. Then for any $f \in V'$, the following boundary value
problem
\begin{eqnarray}
&& u - \Delta u + \nabla p = f \ {\rm in} \ \Omega\\
&& \nabla\cdot u = 0 \ {\rm in} \ \Omega\\
&& u\cdot n = 0, \ 2(S(u)n)_\tau =-\gamma u_\tau + b \ {\rm on} \ \partial\Omega
\end{eqnarray}
has a unique solution $u\in V$ in the sense that
\begin{eqnarray}
(u,\phi) + \int_{\partial\Omega}(\gamma u + b)\cdot \phi +
2\int_\Omega S(u)\cdot S(\phi)= (f,\phi), \forall \phi\ {\rm in}\
V
\end{eqnarray}
Moreover, if $f\in H$ and $b\in H^{\frac{1}{2}}(\partial\Omega)$,
then $u\in H^2(\Omega)$.
\end{thm}

\section{Functional setting of the LNS-$\alpha$ equation}

In this section, we formulate the following
boundary value problem for the LNS-$\alpha$ system:
\begin{eqnarray}
&& \partial_t v - \Delta v + \nabla\times v \times u +
\nabla p=0 \ {\rm in} \ \Omega\\
&& \nabla\cdot v = 0 \ {\rm in} \ \Omega\\\label{4.1}
&& u -
\alpha \Delta u +
\nabla \tilde{p} = v \ {\rm in} \ \Omega\\
&& \nabla\cdot u = 0 \ {\rm in} \ \Omega
\end{eqnarray}
with the VBS conditions
\begin{eqnarray}
&& v \cdot n = 0,\  n\times \nabla\times v = \beta v  \ {\rm on} \ \partial\Omega\\
&& u \cdot n = 0,\  n\times \nabla\times u = \beta u  \ {\rm on}
\ \partial\Omega
\end{eqnarray}
Due to theorem 3.5., $A_{\alpha} = I - \alpha
P\Delta$ is also a positive definite self-adjoint operator with
domain
$\emph{D}(A_\alpha) =W_\beta\cap H$ for any $\alpha > 0$. We have

\begin{prop} The linear operator $T_\alpha =
A_{\alpha}^{-1}: H\mapsto W_\beta\cap H$ is well defined with $u =
T_\alpha v\in W_\beta\cap H$ given by the Stokes boundary value
problem
\begin{eqnarray}
&&  u - \alpha\Delta u +
\nabla \tilde{p} = v \ {\rm in} \ \Omega\\
&& \nabla\cdot u = 0 \ {\rm in} \ \Omega\\
&& u\cdot n = 0,\  n\times(\nabla\times u) = \beta u \
{\rm on} \
\partial\Omega
\end{eqnarray}
and is bounded, i.e.
\begin{eqnarray}
 \|u\|_2\leq c_\alpha \|v\|
\end{eqnarray}
for some constant $c_\alpha$ depending on $\alpha$.
\end{prop}
We now estimate the nonlinearity.
Let $v\in V\subset H$ so that $T_\alpha v$ is defined. Set
\begin{eqnarray}
&& B(v,u) = P(\nabla\times v \times u), \quad \forall\,u\in W_\beta, \quad v\in V\\
&&      B_\alpha(v) = B(v, T_\alpha v), \quad v\in V
    \end{eqnarray}
for $\alpha > 0$. Then we have

\begin{lem} The nonlinearity $B_\alpha (v): V \mapsto
H$ is locally Lipshitz for $\alpha >0$.
\end{lem}
\textbf{Proof:} Clearly, $B_\alpha$ is well-defined due to (4.10). For any $v_1$, $v_2 \in V$,
\begin{equation}
 \|B_\alpha(v_1)-B_\alpha(v_2)\| \leq \|\nabla\times (v_1 -v_2)\times T_\alpha v_1 -
 \nabla\times v_2 \times T_\alpha (v_1-v_2)\|
    \end{equation}
Note that for all $\phi\in V$, $\psi\in L^\infty(\Omega)$,
\begin{equation}
 \|\nabla\times (\phi)\times \psi\|\leq c
 \|\phi\|_1\|\psi\|_{L^\infty(\Omega)}
    \end{equation}
and
\begin{equation}
 \|w\|_{L^\infty(\Omega)}^2 \leq c \|w\|_1\|w\|_2
    \end{equation}
It follows that
\begin{equation}
 \|B_\alpha(v_1)-B_\alpha(v_2)\| \leq c (\|v_1\|_1
 +\|v_2\|_1)\|(v_1-v_2)\|_1
    \end{equation}
which implies the lemma.
\\

We now can formulate the initial boundary problem of the
LNS-$\alpha$ equations (4.1)-(4.6) as an abstract equation
\begin{eqnarray}
&&   v^{\prime} - P\Delta v  + B(v,u) = 0\\
&&  u = T_\alpha v
    \end{eqnarray}
on $H$, with a parameter $\alpha\in (0,\infty)$.
\\

The weak solutions of the initial boundary problem can be defined as below.

\begin{defn} $(v,u)$ is a weak solution of (4.1)-(4.6)
with $\alpha > 0$ for LNS-$\alpha$ equations with initial data $v_{0}\in H$ on
the time interval $[0,T)$ if $v\in L^{2}(0,T; V)\cap
C_{w}([0,T];H)$, $v^{\prime}\in L^{1}(0,T; V')$, $u\in
L^{2}(0,T; V\cap H^3(\Omega))\cap C_{w}([0,T];W_\beta)$, $u^{\prime}\in L^{1}(0,T; V)$ such that
\begin{eqnarray}
&&   (v^{\prime}, w)  +  a_{\beta}(v,w) + (B(v,u),w) = 0,\  a.e. \ t\in [0,T)\\
&&   u = T_\alpha v,\  a.e. \ t\in [0,T)
    \end{eqnarray}
for all $w\in V$.
\end{defn}
For the special case $\alpha=0$, we define also the corresponding weak solutions for the NS as follows

\begin{defn}\ $(v,u)$ is a weak solution of
(4.1)-(4.6) with $\alpha = 0$ (NS equations) initial
data $v_{0}\in H$ on the time interval $[0,T)$ if
$v, u\in L^{2}(0,T; V)\cap C_{w}([0,T];H)$ and
$v^{\prime},v^{\prime}\in L^{1}(0,T; V')$ such that
\begin{eqnarray}
&&   (v^{\prime}, w)  +  a_{\beta}(v,w) + (B(v,u),w) = 0,\  a.e. \ t\in [0,T)\\
&&   u = v,\  a.e. \ t\in [0,T)
    \end{eqnarray}
for all $w\in V$.
\end{defn}
For later use, one can also define the fractional powers of the operator $A_\beta = I-P\Delta$ in
theorem 3.5., $A^s_\beta: D(A^s_\beta)\mapsto H$ for $s\geq 0$ by
\begin{eqnarray}
A^s_\beta\,u = \sum^\infty_{j=1} (1+\lambda_j)^s\,u_j\,e_j
\end{eqnarray}
for $u\in\sum^\infty_{j=1}\,u_j\,e_j\in D(A^s_\beta)$, where
\begin{eqnarray}
D(A^s_\beta)=\left\{ u=\sum^\infty_{j=1}\,u_j\,e_j; \sum^\infty_{j=1} (1+\lambda_j)^{2s}\,|u_j|^2<\infty\right\}
\end{eqnarray}
equipped with the graph norm
$$||u||^2_{D(A^s_\beta)}=(A^s_\beta\,u, A^s_\beta\,u).$$
It can be checked easily that $A^s_\beta: D(A^{s+t}_\beta)\mapsto D(A^t_\beta)$ is an isomorphism for all $s$, $t\geq 0$, $D(A^1_\beta)=D(A_\beta)=H^2(\Omega)\cap W_\beta$, and $D(A^{\frac{1}{2}}_\beta)=V$ with equivalent norms $||u||_{D(A^{\frac{1}{2}}_\beta)}$ and $H^1(\Omega)$-norm. Denote by $D(A^{-s}_\beta)$ the dual space of $D(A^s_\beta)$ for any $s\geq 0$. Then the operator $A^s_\beta$ can be extended to an operator: $H\mapsto D(A^{-s}_\beta)$ by
\begin{eqnarray}
(A^s_\beta\,u,v)=(u,A^s_\beta\,v), \qquad \forall u\in H, \quad v\in D(A^s_\beta)
\end{eqnarray}
It follows from the definition that
\begin{eqnarray}
||u||^2_{D(A^{-s}_\beta)}=||A^{-s}_\beta\,u||^2, \qquad \forall u\in D(A^{-s}_\beta)
\end{eqnarray}
and $A^s: D(A^{s+t}_\beta)\mapsto D(A^t_\beta)$ is an isomorphism for $s$, $t\in\mathbb{R}$, and furthermore,
\begin{eqnarray}
||A^{\frac{s+t}{2}}_\beta\,u||^2 = (A^s_\beta\,u,A^t_\beta\,u)\leq ||A^s_\beta\,u||\,||A^t_\beta\,u||, \, \forall u\in D(A^s_\beta)\cap D(A^t_\beta)
\end{eqnarray}
holds true for all $s$, $t\in\mathbb{R}$.

\section{Well-Posedness of the LNS-$\alpha$ Equations}

In this section, we investigate the well-posedness of the initial boundary value problem of
the LNS-$\alpha$ equations (4.1)-(4.6) by a Gelerkin approximation based on the
orthogonal basis given in theorem 3.5..

\subsection{Local well-posedness}

We start with the following local well-posedness result.

\begin{thm} Let $v_{0}\in H$ and $\alpha > 0$. Then there is a
time $T^{*} = T^{*}(v_{0})> 0$ such that the problem (4.1)-(4.8) has a unique
weak solution of $(v,u)$ with initial data $v_0$ on the interval $[0,T^{*})$ in the sense
of definition 4.1 for any $T\in(0,T^{*})$, which satisfies the energy
equation
\begin{equation} \label{a19}
  \frac{d}{dt}\|v\|^{2} + 2 a_{\beta}(v,v) + (B(v,u),v) = 0,\
  {\rm on}\ [0,T]
    \end{equation}
in the sense of distribution. Furthermore, if $v_{0}\in V$, then
\begin{eqnarray}
&& v\in L^{2}(0, T ;W_{\beta}\cap H)\cap C([0, T) ; V) \\
&& v^{\prime}\in L^{2}(0, T ; H)
\end{eqnarray}
and the energy
equation
\begin{equation}
  \frac{d}{dt}a_{\beta}(v,v) + 2 \|P\Delta v\|^{2}
  + 2(B(v),-\Delta v) = 0
    \end{equation}
is valid.
\end{thm}
\textbf{Proof:} Let $ v_{0} \in H $. Consider the following system of
ordinary differential equations
\begin{eqnarray}
&&  v_{j}^{\prime}(t) +  \lambda_{j}v_{j}(t) + g_{j}(\mathcal{V})= 0\\
&& v_{j}(0) = (u_{0},e_{j})
    \end{eqnarray}
$j = 1,\cdots m$, where $\mathcal{V} = (v_{j})$ and
\begin{eqnarray}
&&  g_{i}(\mathcal{V}) = (B(\Sigma_{1}^{m} v_{j}e_{j}, u_m),e_{i}) \\
&&  u_m = T_\alpha(\Sigma_{1}^{m} v_{j}e_{j})
    \end{eqnarray}
Note that all norms are equivalent in a finite dimensional linear space. It follows
from lemma 4.2. that $(g_{j}(\mathcal{V}))$ is locally Lipshitz
in $\mathcal{V}$ and thus the systems is locally well posed and
equivalent to the following partial differential equations
\begin{eqnarray}
&&\label{a9}  v_{m}^{\prime}(t,x) - P\Delta v_{m}(t,x) + P_{m}B (v_{m},u_m)(t,x)= 0 \\
&& u_m = T_\alpha(\Sigma_{1}^{m} v_{j}e_{j})\\
&&\label{a10} v_{m}(0) = P_{m}(v_{0})
    \end{eqnarray}
where $v_{m}(t,x)=\Sigma_{1}^{m} v_{j}(t)e_{j}(x)$, and
$P_{m}$ is the orthogonal projection of $H$ onto the space
$spin\{e_{j}\}_{1}^{m}$.\\

Taking the inner product of (5.9) with $v_{m}$ and noting that
\begin{equation}
 (P_{m}B(v_{m},u_{m}),v_{m}) =
 \int_{\Omega}\nabla\times v_{m}\times u_{m}\cdot(v_{m})dx
    \end{equation}
one can get
\begin{equation}
  \frac{d}{dt}\|v_{m}\|^{2} + 2a_\beta(v_m,v_m) + (B(v_m,u_m),v_m) = 0
    \end{equation}
It follows from the definition of $T_\alpha$ that
\begin{equation}
  |(B(v_m,u_m),v_m)|\leq
  c\|v_m\|_1\|v_m\|\|u_m\|_{L^\infty(\Omega)}\leq c
  \|v_m\|_1\|v_m\|^2
    \end{equation}
Note that
\begin{equation}
  \|\phi\|^2\leq c\|\phi\|_1^2\leq c(\|\phi\|^2+
  a_{\beta}(\phi,\phi))
    \end{equation}
for all $\phi\in V$. It follows that
\begin{equation} \label{5.3}
  \frac{d}{dt}\|v_{m}\|^{2} + a_{\beta}(v_m,v_m) \leq c(\|v_{m}\|^{2} +1)\|v_{m}\|^{2}
    \end{equation}
Hence, there is a time $T>0$ such that
\[ \{v_{m}\}\ \ is\ bounded\ in\   L^{\infty}(0, T;
   H) \]
\[ \{v_{m}\}\ \ is\ bounded\ in\   L^{2}(0, T;
  V) \]
Note that for $\phi \in V $,
\begin{equation} \label{a20}
   |(A_\beta v_{m},\phi)|  \leq |(v_m, \phi)| + |a_{\beta}(v_m,\phi)|
    \end{equation}
which implies that
\begin{equation}
 \{A_\beta v_{m}\}\ is\ bounded\ in\ \L^{2}(0,T; V')
  \end{equation}
Since
\begin{equation}
 \|u_m\|_{L^\infty(\Omega)}= c\|T_\alpha v_m\|_2\leq c\|v_m\|
  \end{equation}
it follows that
\begin{equation}
   |(P_{m}B(v_m,u_{m}),\phi)|  = |(\nabla\times v_{m}\times u_{m},P_m\,\phi)|
   \leq C\|v_m\|_1\|v_m\|\|\phi\|_1
    \end{equation}
for all $\phi\in V$, which implies that
\begin{equation}
  \{P_m\,B(v_{m}, u_m)\} \ is\ bounded\ in\ \L^{2}(0,T; V')
  \end{equation}
Hence
\begin{equation}
  \{v_{m}^{\prime}\} \ is\ bounded\ in\ \L^{2}(0,T; V')
  \end{equation}
By using a similar argument in \cite{C88}, it shows that there
is a subsequence also denoted by $v_m$ and a
$v\in  L^{\infty}(0, T; H)\cap L^{2}(0, T; V)$ such that
\begin{eqnarray}
&& v_{m}\rightarrow v  \ in\ L^{\infty}(0, T ;H)\ weak-star\\
&& v_{m}\rightarrow v \ in\ L^{2}(0, T ; V)\ weakly\\
&& v_{m}\rightarrow v \ in\ L^{2}(0, T ;H)\ strongly
    \end{eqnarray}
Consequently, $u_m = T_\alpha v_m$ has the property:
\begin{eqnarray}
&& u_{m}\rightarrow T_\alpha v  \ in\ L^{\infty}(0, T ;W_\beta\cap H)\ weak-star\\
&& u_{m}\rightarrow T_\alpha v \ in\ L^{2}(0, T ; V\cap H^3(\Omega))\ weakly\\
&& u_{m}\rightarrow T_\alpha v \ in\ L^{2}(0, T ;W_\beta\cap H)\
strongly
    \end{eqnarray}
Passing to the limit of a subsequence, it is showed that $(v,u)$ is a weak
solution in the sense of definition 4.3..
It also follows that the energy equation
\begin{equation} \label{5.1}
  \frac{d}{dt}\|v\|^{2} + 2 a_{\beta}(v,v) + (B_\alpha(v),v) = 0
    \end{equation}
is valid on the interval $[0,T]$ in the sense of distribution.
\\

Let $v_1$ and $v_2$ be any two solutions. Then $w =
v_1-v_2$ satisfies the following equation
\begin{eqnarray}
&&\label{a9}  w ^{\prime}- P\Delta w + P(B_\alpha(v_{1})-B_\alpha(v_2))= 0 \\
&&\label{a10}w(0) = 0
    \end{eqnarray}
and the energy equation
\begin{equation}
  \frac{d}{dt}\|w\|^{2} + 2 a_{\beta}(w,w) + (B_\alpha(v_1)-B_\alpha(v_2),w) = 0
    \end{equation}
It follows from the local Lipshitz continuity stated in lemma 4.2.
and the Gronwall inquality that
\begin{equation}
  \|w\|^{2} \leq c(T)\|w(0)\|^2 \quad {\rm on} [0,T]
    \end{equation}
which implies the uniqueness of the solution. Consequently, the
convergence of the whole sequence follows.\\

By the standard continuation method, there is a $T^*> 0$
such that the weak solution does exist on $[0,T]$ for all $T < T^*$,
and if $T^* < \infty$ then
\[ \|v(t)\|\rightarrow \infty, as \ t\rightarrow T^*\]
Let $v_0\in V$. Taking the inner product of (5.9) with $-P\Delta v_{m}$ and noting that
\begin{equation}
 (P_{m}B(v_{m},u_m),-P\Delta v_{m}) = (B(v_{m},u_m),-P\Delta v_{m})
    \end{equation}
one gets
\begin{equation}
  \frac{d}{dt}a_{\beta}(v_m,v_m) + 2 \|-P\Delta v_m\|^2 + (B_\alpha(v_m),-P\Delta v_m) = 0
    \end{equation}
It follows that
\begin{equation}
  \frac{d}{dt}a_{\beta}(v_m,v_m) +  \|-P\Delta v_m\|^2
  \leq \|v_m\|_1^2\|u_m\|_{L^\infty(\Omega)}^2
    \end{equation}
Due to
\begin{equation}
  \|u_m\|_{L^\infty(\Omega)}\leq c_\alpha\|v_m\|
    \end{equation}
\\
$v_0\in V$, the bounds of $v_m$ in $L^2(0,T; V)$, and
the Gronwall's inequality, one has
\[ \{v_{m}\}\ \ is\ bounded\ in\   L^{\infty}(0, T;
   V)\]
\[ \{v_{m}\}\ \ is\ bounded\ in\   L^{2}(0, T;
  W_{\beta}\cap H) \]
which, together with the uniqueness, implies that the whole sequence
indeed converges in the sense
\begin{eqnarray}
&& v_{m}\rightarrow v  \ in\ L^{\infty}(0, T ;V)\ weak-star\\
&& v_{m}\rightarrow v \ in\ L^{2}(0, T ; W_\beta\cap H)\ weakly\\
&& v_{m}\rightarrow v \ in\ L^{2}(0, T ;V)\ strongly
    \end{eqnarray}
This completes the proof of theorem 5.1..

\subsection{Global well-posedness}

Now, we prove the following global well-posedness result.

\begin{thm} If $v_{0}\in H, \alpha > 0$, then the
solution $v$ obtained in theorem 5.1. is global, i.e., $T^{*} =
T^{*}(v_{0})=\infty$.
\end{thm}
\textbf{Proof:} Let $v$ be the weak solution on the interval
$[0,T]$. Then,
\begin{eqnarray}
&& u= T_\alpha v\in L^{\infty}(0, T ;W_\beta\cap H)
\end{eqnarray}
Taking $u$ as a test function yields
\begin{equation}
  (v^{\prime},u)  +  a_{\beta}(v,u) + (B(v,u),u) = 0
    \end{equation}
Since $v = (I -\alpha P\Delta) u$, then
\begin{equation}
  2 (v^{\prime},u)  = \frac{d}{dt}(\|u\|^2 + \alpha a_\beta(u,u))
    \end{equation}
in the sense of distribution on $[0,T]$. Note that
\begin{equation}
  (B(v,u),u) = \int_\Omega (\nabla\times v)\times u\cdot u = 0
    \end{equation}
It follows that
\begin{equation}
\frac{d}{dt}(\|u\|^2 + \alpha a_\beta(u,u))+
2(\int_{\partial\Omega}\beta u\cdot v + \int_\Omega (\nabla\times
v)\cdot (\nabla\times u)) = 0
    \end{equation}
Due to the smoothness and the boundary condition for $u$, it holds that
\begin{equation}
\int_\Omega (\nabla\times v)\cdot (\nabla\times u) =
-\int_{\partial\Omega}\beta u\cdot v + \int_\Omega (-\Delta
u)\cdot v
    \end{equation}
Consequently
\begin{equation} \label{5.2}
\frac{d}{dt}(\|u\|^2 + \alpha a_{\beta}(u,u))+ 2(a_\beta(u,u)+
\alpha\|P\Delta u\|^2) = 0
    \end{equation}
It follows that
\begin{equation}
(\|u\|^2 + \alpha a_{\beta}(u,u))\leq (\|u_0\|^2 + \alpha
a_{\beta}(u_0,u_0))
    \end{equation}
and
\begin{equation}
\int_0^t (a_\beta(u,u)+ \alpha\|P\Delta u\|^2)d\tau\leq
(\|u_0\|^2 + \alpha a_{\beta}(u_0,u_0))
    \end{equation}
On the other hand, it follows from the energy equation (5.1) and a similar argument as for (5.17) that
\begin{equation}
  \frac{d}{dt}\|v\|^{2} + a_{\beta}(v,v) \leq c \|v\|^{4} + 1
    \end{equation}
Noting that
\begin{equation}
  \|v\|^{2}  \leq c (\|u\|^2 + \alpha^2\|P\Delta u\|^2)
    \end{equation}
it follows that
\begin{equation}
\|v\|^2 + \int_0^t a_{\beta}(v,v)\leq c
    \end{equation}
for some constant $c$ depending only on $v_0$ and $\alpha$. Thus $T^* =
\infty$. The theorem is proved.

\section{Vanishing $\alpha$ Limit and the NS Equations}

In this section, we investigate the vanishing $\alpha$ limit of
the solutions of the LNS-$\alpha$ equations ($\alpha\rightarrow
0$) to that of the NS equations. We will prove both weak
and strong convergence results. Then, the global
existence of weak solutions and the local unique strong solution
to the NS equations with the VSB condition are followed.

\subsection{Weak Convergence and Global Weak Solutions of the NS}

We first prove

\begin{thm} Let $v_{0}\in H$, and
$(v^\alpha,u^\alpha)$ be the global weak solution stated in theorem 5.2.
corresponding to the parameter $\alpha>0$. Then for any given $T>0$
there is a subsequence $u^{\alpha_j}$ of $u^\alpha$ and a $(v^0, u^0)$
satisfying
\begin{eqnarray}
&& v^0 \in L^{2}(0, T ;V)\cap C_w([0, T] ; H) \\
&& (v^0)^{\prime}\in L^{\frac{4}{3}}(0, T ; V')
\end{eqnarray}
such that
\begin{eqnarray}
&& v^{\alpha_j}\rightarrow v^0 \ in\ L^{2}(0, T ; H)\ weakly\\
&& v^{\alpha_j}\rightarrow v^0 \ in\ L^{2}(0, T
;D(A_\beta^{-\frac{1}{4}}))\ strongly
    \end{eqnarray}
\begin{eqnarray}
&& u^{\alpha_j}\rightarrow v^0 \ in\ L^{2}(0, T ; V_\beta)\ weakly\\
&& u^{\alpha_j}\rightarrow v^0 \ in\ L^{2}(0, T
;D(A_\beta^{-\frac{1}{4}}))\ strongly
    \end{eqnarray}
Moreover $(v^0,v^0)$ is a weak solution of the initial boundary problem
of the NS equations (4.1)-(4.6) with $\alpha =0$ and satisfies the energy
inequality
\begin{equation}
\frac{d}{dt}\|v^0\|^2 + 2 a_\beta(v^0,v^0)\leq 0
    \end{equation}
\end{thm}
\textbf{Proof:} Let $v_0\in H$, $T>0$, and $(v^\alpha,u^\alpha)$ be the global weak solution to (4.1)-(4.6)
corresponding to $1\geq\alpha>0$. It follows from (\ref{5.2}) that
\begin{equation} \label{6.1}
 \|u^\alpha \|^2 + \alpha a_\beta(u^\alpha,u^\alpha)+
 \int_0^t(a_\beta(u^\alpha,u^\alpha) + \alpha\|P\Delta u^\alpha \|^2)d\tau\leq c
    \end{equation}
for some constant $c$ independent of $\alpha$.
For any $\phi\in W_\beta\cap H$, we have
\begin{equation} \label{6.5}
  (B(v^\alpha,u^\alpha),\phi) = \int_{\Omega}(\nabla\times v^\alpha\times
  u^\alpha)\phi dx = I +II
    \end{equation}
where
\begin{equation}
  I = \int_{\partial\Omega}(n\times v^\alpha)\cdot(u^\alpha\times
  \phi) dS
    \end{equation}
\begin{equation}
  II = \int_{\Omega}v^\alpha\cdot(-
  u^\alpha\cdot\nabla \phi-\phi\cdot\nabla u^\alpha) dx
    \end{equation}
Since $u\cdot n = 0$ and $\phi\cdot n = 0$ on the boundary so
\begin{equation}
  u^\alpha\times \phi = \lambda n \ {\rm on}\ \partial\Omega
    \end{equation}
Hence
\begin{equation} \label{6.6}
  I = 0
    \end{equation}
To estimate $II$, we note that
\begin{equation}\label{6.1.2}
  |\int_{\Omega}v^\alpha\cdot(
  u^\alpha\cdot\nabla \phi)dx|\leq c(\|u^\alpha\|+\alpha\|P\Delta
  u^\alpha\|)\|u^\alpha\|_{L^3(\Omega)}\|\nabla \phi\|_{L^6(\Omega)}
    \end{equation}
\begin{equation}\label{6.1.3}
   \|u^\alpha\|_{L^3(\Omega)}^2\leq
   c\|u^\alpha\|\|u^\alpha\|_1\leq c\|u^\alpha\|^\frac{3}{2}(\|u^\alpha\|+\|P\Delta
   u^\alpha\|)^\frac{1}{2}
    \end{equation}
\begin{equation} \label{6.1.4}
  \|\nabla \phi\|_{L^6(\Omega)}\leq c  \|A_\beta \phi\|
    \end{equation}
Then, due to (6.8), it holds that
\begin{equation}
  |\int_{\Omega}v^\alpha\cdot(
  u^\alpha\cdot\nabla \phi)dx|\leq c
  ((a_\beta(u^\alpha,u^\alpha))^\frac{1}{2}+ \alpha||P\Delta u^\alpha|| + \alpha||P\Delta
  u^\alpha\|^\frac{5}{4})\|A_\beta \phi\|
    \end{equation}
Next,
\begin{equation}
  |\int_{\Omega}v^\alpha\cdot(
  \phi\cdot\nabla u^\alpha)dx|\leq c(\|u^\alpha\|+\alpha\|P\Delta
  u^\alpha\|)\|u^\alpha\|_1\| \phi\|_{L^\infty(\Omega)}
    \end{equation}
which implies that
\begin{equation}
  |\int_{\Omega}v^\alpha\cdot(
  \phi\cdot\nabla u^\alpha)dx|\leq c
  ((a_\beta(u^\alpha,u^\alpha))^\frac{1}{2}+ \alpha||P\Delta u^\alpha|| + \alpha||P\Delta
  u^\alpha\|^\frac{3}{2})\|A_\beta \phi\|
    \end{equation}
Then for $\alpha< 1$,
\begin{equation} \label{6.2}
 |(B(v^\alpha,u^\alpha),\phi)|\leq c
  (1 + (a_\beta(u^\alpha,u^\alpha))^\frac{1}{2}+ \alpha^\frac{3}{4}\|P\Delta
  u^\alpha\|^\frac{3}{2})\|A_\beta \phi\|
    \end{equation}
It follows from (\ref{6.1}) and (\ref{6.2}) that
$B(v^\alpha,u^\alpha)$ and then $\frac{d}{dt}(v^\alpha)$ are
uniformly bounded in $L^{\frac{4}{3}}(0,T; D(A_\beta^{-1}))$. It follows from (4.7)-(4.9) that
$$(1-\alpha)u^\alpha_t+\alpha A_\beta(u^\alpha_t)=v^\alpha_t,$$
which yields immediately
$$(1-\alpha)||A^{-1}_\beta\,u^\alpha_t||^2 + \alpha||A^{-1}_\beta\,u^\alpha_t||^2 = ||A^{-1}_\beta\,v^\alpha_t||^2.$$
Then
$$||A^{-1}_\beta\,u^\alpha_t||^2 \leq 2||A^{-1}_\beta\,v^\alpha_t||^2$$
for $0<\alpha\leq\frac{1}{2}$. This shows that
$\partial_t\,u^\alpha$ are uniformly bounded in
$L^{\frac{4}{3}}(0,T; D(A^{-1}_\beta))$ as $\partial_t\,v^\alpha$ are. Note
that (\ref{6.1}) also implies that $(u^\alpha)$ are uniformly bounded
in $L^{2}(0,T; V)$ and the duality between $V =
D(A_\beta^{\frac{1}{2}})$ and $D(A_\beta^{-1})$ with respect to the inner
product of $D(A_\beta^{-\frac{1}{4}})$, i.e.,
\[(A_\beta^{-\frac{1}{4}}u,A_\beta^{-\frac{1}{4}}\phi) = (A_\beta^{\frac{1}{2}}u, A_\beta^{-1}\phi)\]
By using the
standard compactness argument (see \cite{FHT,C88}), one can show
that there exist a subsequence $u^{\alpha_j}$ of $u^\alpha$ and a
$v^0$ such that
\begin{eqnarray} \label{6.3}
&& u^{\alpha_j}\rightarrow v^0 \ in\ L^{2}(0, T ; V)\ weakly\\
\label{6.4} && u^{\alpha_j}\rightarrow v^0 \ in\ L^{2}(0, T
;D(A_\beta^{-\frac{1}{4}}))\ strongly
    \end{eqnarray}
Note that
\begin{equation} \label{6.11}
 |(B(v^\alpha,u^\alpha)-B(v^0,v^0),\phi)|\leq I + II
    \end{equation}
where
\begin{equation}
  I = |(B(u^\alpha-v^0,u^\alpha)+ B(v^0,u^\alpha-v^0),\phi)|
    \end{equation}
\begin{equation}
 II = \alpha|(B(P\Delta u^\alpha ,u^\alpha),\phi)|
    \end{equation}
Similar to (\ref{6.5}) and (\ref{6.6}), integrating by part yields
\begin{equation}
|(B(u^\alpha-v^0,u^\alpha),\phi)|
=|\int_\Omega (u^\alpha-v^0)\cdot(u^\alpha\cdot\nabla\phi+\phi\cdot\nabla u^\alpha)|
    \end{equation}
Note that
\begin{equation}
|\int_\Omega (u^\alpha-v^0)\cdot(u^\alpha\cdot\nabla\phi)|\leq c
 \|u^\alpha-v^0\|\|u^\alpha\|^{\frac{1}{2}}\|u^\alpha\|_{L^6(\Omega)}^{\frac{1}{2}}\|\nabla\phi\|_{L^6(\Omega)}
    \end{equation}
and
\begin{equation}
 \|u^\alpha-v^0\|^2\leq c \|u^\alpha-v^0\|_{D(A_\beta^{-\frac{1}{4}})}\|u^\alpha-v^0\|_{D(A_\beta^{\frac{1}{4}})}
    \end{equation}
\begin{equation}
 \|u^\alpha-v^0\|_{D(A_\beta^{\frac{1}{4}})}^2\leq c \|u^\alpha-v^0\|\|u^\alpha-v^0\|_1
    \end{equation}
This, together with (\ref{6.1}), shows that
\begin{equation}
|\int_\Omega (u^\alpha-v^0)|\cdot(u^\alpha\cdot\nabla\phi)\leq
 c\|u^\alpha-v^0\|_{D(A_\beta^{-\frac{1}{4}})}^\frac{1}{2}
 \|u^\alpha-v_0\|_1^{\frac{1}{4}}\|u^\alpha\|_1^{\frac{1}{2}}\|A_\beta\phi\|
    \end{equation}
Hence
\begin{equation}
|\int_\Omega (u^\alpha-v^0)|\cdot(u^\alpha\cdot\nabla\phi)\leq
 c\|u^\alpha-v^0\|_{D(A_\beta^{-\frac{1}{4}})}^\frac{1}{2}
 (\|u^\alpha \|_1^{\frac{3}{4}}+\|v^0\|_1^{\frac{3}{4}})\|A_\beta\phi\|
    \end{equation}
While
\begin{equation}
|\int_\Omega (u^\alpha-v^0)\cdot(\phi\cdot\nabla u^\alpha)|\leq
 \|u^\alpha-v^0\|\|u^\alpha\|_1 \| \phi\|_{L^\infty(\Omega)}
    \end{equation}
It follows that
\begin{equation}
|\int_\Omega (u^\alpha-v^0)|\cdot(u^\alpha\cdot\nabla\phi)\leq
 c\|u^\alpha-v^0\|_{D(A_\beta^{-\frac{1}{4}})}^\frac{1}{2}
 (\|u^\alpha \|_1^{\frac{5}{4}}+ \|v^0\|_1^{\frac{5}{4}})\|A_\beta\phi\|
    \end{equation}
Then
\begin{equation} \label{6.12}
|(B(u^\alpha-v^0,u^\alpha),\phi)|\leq c\|u^\alpha-v^0\|_{D(A_\beta^{-\frac{1}{4}})}^\frac{1}{2}
 (1+ \|u^\alpha \|_1^{\frac{5}{4}}+ \|v^0\|_1^{\frac{5}{4}})\|A_\beta\phi\|
    \end{equation}
Similarly, one can obtain
\begin{equation}
|(B(v^0,u^\alpha-v^0),\phi)|\leq c\|u^\alpha-v^0\|_{D(A_\beta^{-\frac{1}{4}})}^\frac{1}{2}
 (1+ \|u^\alpha \|_1^{\frac{5}{4}}+ \|v^0\|_1^{\frac{5}{4}})\|A_\beta\phi\|
    \end{equation}
It follows that
\begin{equation} \label{6.12}
I\leq c\|u^\alpha-v^0\|_{D(A_\beta^{-\frac{1}{4}})}^\frac{1}{2}
 (1+ \|u^\alpha \|_1^{\frac{5}{4}}+ \|v^0\|_1^{\frac{5}{4}})\|A_\beta\phi\|
    \end{equation}
Similarly,
\begin{equation}
|(B(P\Delta u^\alpha,u^\alpha),\phi)| =|\int_\Omega (P\Delta
u^\alpha)\cdot(u^\alpha\cdot\nabla\phi+\phi\cdot\nabla u^\alpha)|
    \end{equation}
Then
\begin{equation}
|(B(P\Delta u^\alpha,u^\alpha),\phi)| \leq c\|P\Delta
u^\alpha\|\|u^\alpha\|_1\|A_\beta\phi\|
    \end{equation}
It follows that
\begin{equation} \label{6.13}
II \leq c\alpha^{\frac{1}{2}}(\alpha\|P\Delta u^\alpha\|^2 +
\|u^\alpha\|_1^2)\|\phi\|_2
    \end{equation}
It follows from (\ref{6.11}),(\ref{6.12}),(\ref{6.13}),
(\ref{6.1}) and (\ref{6.3}) that
\begin{eqnarray}
&& B(v^{\alpha_j},u^{\alpha_j}) \rightarrow B(v^0,v^0) \ in\
L^{1}(0, T ; D(A_\beta^{-1}))\ strongly
    \end{eqnarray}
which enables us to pass the limit in (4.19)-(4.20) to show that $v^0$ satisfies
\begin{equation}\label{6.1.5}
   ((v^0)^{\prime}, \phi)  + a_{\beta}(v^0,\phi) +
((\nabla\times v^0)\times v^0,\phi) = 0,\  a.e. \ t
    \end{equation}
for all $\phi\in C^\infty(\Omega)\cap V$ in the sense of
distribution on $[0,T]$. Note that $v^0\in L^{2}(0, T ; V)$
implies
$(v^0)^\prime \in L^{\frac{4}{3}}(0, T ; V')$. Thus (\ref{6.1.5}) is also valid
for all $\phi\in V$.
\\
Due to (\ref{5.2}), it holds that
\begin{equation}
\frac{d}{dt}(\|u^\alpha\|^2 + \alpha
a_{\beta}(u^\alpha,u^\alpha))+ 2 a_\beta(u^\alpha,u^\alpha) \leq
0
    \end{equation}
Passing to the limit and noting the weak lower
semi-continuity of the norm, one gets
\begin{equation}
\frac{d}{dt}\|v^0\|^2 + 2 a_\beta(v^0,v^0)\leq 0
    \end{equation}
Note that
\begin{eqnarray}
&&   (v^\alpha-v^0, \phi)= (u^\alpha-v^0, \phi) +
\alpha((A_\beta^{\frac{3}{4}}u^\alpha,A_\beta^{\frac{1}{4}}\phi)-(u^\alpha,\phi))
    \end{eqnarray}
for $\phi\in D(A_\beta^{-\frac{1}{4}})$. Then
\[  \|v^\alpha-v^0\|_{D(A_\beta^{-\frac{1}{4}})}^2\leq
\|v^\alpha-v^0\|_{D(A_\beta^{-\frac{1}{4}})}^2 +
c\alpha^{\frac{1}{2}}(\alpha\|P\Delta u^\alpha\|^2 +
\|u^\alpha\|_1^2)
   \]
It follows that
\begin{eqnarray}
&& v^{\alpha_j}\rightarrow v^0 \ in\ L^{2}(0, T
;D(A_\beta^{-\frac{1}{4}}))\ strongly
    \end{eqnarray}
Note that
\begin{eqnarray}
&&   (v^\alpha - v^0, \phi) = (u^\alpha-v^0, \phi) -
\alpha(P\Delta u^\alpha, \phi)
    \end{eqnarray}
It follows from (\ref{6.1}) and (\ref{6.3}) that
\begin{eqnarray}
&& v^{\alpha_j}\rightarrow v^0 \ in\ L^{2}(0, T ; H)\ weakly
    \end{eqnarray}
Hence, the theorem is proved.

\subsection{Strong Convergence and the Strong Solutions of the NS}

We now turn to the strong convergence of the strong solutions of
the LNS-$\alpha$ to that of the NS equations, and prove

\begin{thm} Let $v_{0}\in V$ and $(v^\alpha,u^\alpha)$
be the strong solution stated in theorem 5.1. corresponding to the
parameter $\alpha>0$.  Then there is a $T>0$
and a $v^0 \ in\ \ L^\infty(0, T ;V)\cap L^{2}(0, T ; W_\beta\cap
H)$ such that
\begin{eqnarray}
&& v^{\alpha}\rightarrow v^0 \ in\ L^{2}(0, T ; W_\beta\cap
H)\ weakly\\
&& v^{\alpha}\rightarrow v^0 \ in\ L^{2}(0, T ;V)\ strongly
    \end{eqnarray}
\begin{eqnarray}
&& u^{\alpha}\rightarrow v^0 \ in\ L^{2}(0, T ; W_\beta\cap
H)\ weakly\\
&& u^{\alpha}\rightarrow v^0 \ in\ L^{2}(0, T ;V)\ strongly
    \end{eqnarray}
with $v^0$ being a weak solution to the initial boundary problem of
the NS  equation (4.1)-(4.6) with $\alpha =0$ which is unique and
thus called the strong solution. Consequently, it can be extended to the maximal
existence time interval $[0,T^*)$ such that if $T^* < \infty$ then
\[\|v^0\|_1 \rightarrow \infty, \ as \ t\rightarrow T^*\]
Moreover, the following energy equation holds:
\begin{equation}
  \frac{d}{dt}a_{\beta}(v^0,v^0) + 2 \|P\Delta v^0\|^{2}
  - 2(B(v^0,v^0),P\Delta v^0) = 0
    \end{equation}
\end{thm}
\textbf{Proof:} It follows from the energy equation (5.5) that
\begin{equation}
  \frac{d}{dt}a_{\beta}(v^\alpha,v^\alpha) + \|P\Delta v^\alpha\|^{2}
  \leq c\|B(v^\alpha,u^\alpha)\|^2
    \end{equation}
Note that
\begin{equation}
  \|B(v^\alpha,u^\alpha)\|^2\leq c\int_{\Omega}|\nabla\times v^\alpha|^2|u^\alpha|^2dx
  \leq c \|\nabla\times
  v^\alpha\|_{L^3(\Omega)}^2\|u^\alpha\|_{L^6(\Omega)}^2
    \end{equation}
\begin{equation}
\|\nabla\times v^\alpha\|_{L^3(\Omega)}^2
\leq c(\|v^\alpha\| + \|P\Delta v^\alpha\|)\|v^\alpha\|_1
    \end{equation}
\begin{equation}
  \|u^\alpha\|_{L^6(\Omega)}\leq c\|u^\alpha\|_1
    \end{equation}
and
\begin{equation}\label{6.2.0}
||u^\alpha||_1 \leq c||v^\alpha||_1
\end{equation}
which follows from the fact that
\begin{eqnarray}
&&\|u^\alpha\|^2 + \alpha a_\beta(u^\alpha,u^\alpha) = (v^\alpha,
u^\alpha)\\
&&  a_\beta(u^\alpha,u^\alpha) + \alpha\|P\Delta u^\alpha\|^{2}
  = a_\beta(u^\alpha,v^\alpha)
    \end{eqnarray}
Consequently,
\begin{equation}\label{6.2.1}
  \frac{d}{dt}a_{\beta}(v^\alpha,v^\alpha) + \frac{1}{2}\|P\Delta v^\alpha\|^{2}
  \leq c(1+\|v^\alpha\|_1^2)||v^\alpha||^4_1
    \end{equation}
Combining this with similar estimates for (5.17) yields
\begin{equation}
  \frac{d}{dt}(\|v^\alpha\|^2+ a_{\beta}(v^\alpha,v^\alpha))
   + \frac{1}{2}(\|v^\alpha\|^2 + \|P\Delta v^\alpha\|^{2})
  \leq c(1+ \tilde{a}_{\beta}(v^\alpha,v^\alpha))^3
    \end{equation}
Comparing it with the following ordinary differential equation
\begin{equation*}
  \frac{d}{dt}y = c(1+y)^3
    \end{equation*}
with $y(0)= \tilde{a}_{\beta}(v_0,v_0)$ shows that there is a time
$T$ such that
\[ v^\alpha\ is\ uniform\ bounded\ in\
L^\infty(0,T; V)\cap L^2(0,T; W_\beta\cap H)\]
It follows from this, (6.55)-(6.60), and (4.19) that
\[(v^\alpha)^\prime\ is\ uniform\ bounded\ in\ L^2(0,T;H).\]
Hence, by using the standard compactness argument, we find a subsequence
$v^{\alpha_j}$ of $v^\alpha$ and a $v^0$ such that
\begin{eqnarray}
&& v^{\alpha_j}\rightarrow v^0 \ in\ L^{2}(0, T ; W_\beta\cap H)\ weakly\\
&& v^{\alpha_j}\rightarrow v^0 \ in\ L^{2}(0, T ;V)\
strongly
    \end{eqnarray}
which enables one to pass to the limit to find $v^0\in C([0, T]
;V)\cap L^{2}(0, T ; W_\beta\cap H)$ such that $(v^0,v^0)$is a
(strong) solution of the NS equations.
\\
Let $v^0_1$ and $v^0_2$ be two strong solutions to the Navier-Stokes equations with same initial data. Set $w =
v^0_1- v^0_2$. Then
\begin{equation}
\frac{d}{dt}\|w\|^2 + 2  a_\beta(w,w) + 2 (B(v^0_1,v^0_1) -
B_0(v^0_2,v^0_2), w) = 0
    \end{equation}
Note that
\begin{equation*}
\begin{array}{rl}
|(B(v^0_1,v^0_1)-B(v^0_2,v^0_2),w)|\leq & |(B(w,v^0_1),w)|+|(B(v^0_2,w),w)|\\
\leq & \tilde{a}_\beta(w,w)+c(||v^0_1||_{L_\infty}(t)+||\nabla\times v^0_2||^4)||w||^2
\end{array}
\end{equation*}
which, together with (5.53) and Gronwall's inequality, yields $||w||=0$. Thus we
have obtained the uniqueness of the strong solution to the
initial boundary value problem for the Navier-Stokes equations. By the standard
continuation method, the strong solution can be extended to the maximum existent
time interval $[0, T^*)\supset [0,T]$, and the energy equation
follows from the
smoothness of the solution. Consequently, the convergence of the whole sequence of $v^\alpha$
follows.
\\

Finally, we prove the convergence of $u^\alpha$. It follows from (\ref{4.1}) that
\begin{eqnarray}
\nabla\times u^\alpha - \alpha \Delta (\nabla\times u^\alpha)
 = \nabla\times v^\alpha, \ {\rm in} \ \Omega
\end{eqnarray}
Taking the inner product of above equality with $-\Delta(\nabla\times u^\alpha)$ and
integrating by part, we can get
\begin{eqnarray}
\|\Delta u^\alpha\|^2+ \alpha \|(\nabla\times)^3 u^\alpha\|
 = (\Delta u^\alpha, \Delta v^\alpha)
 + \int_{\partial\Omega}\Delta u^\alpha\cdot
 \beta(v^\alpha-u^\alpha)
\end{eqnarray}
To handle the last term on the right hand side above, we use the fact
$v^\alpha-u^\alpha=n\times((v^\alpha-u^\alpha)\times n)$ on $\partial\Omega$
and the Stokes formula to get
\begin{eqnarray}
\begin{array}{cl}
& \displaystyle \int_{\partial\Omega} \Delta u^\alpha\cdot(\beta(v^\alpha-u^\alpha))
= \int_{\partial\Omega} \Delta u^\alpha\cdot(n\times(\beta(v^\alpha-u^\alpha)\times n))\\[3mm]
= & \displaystyle \int_{\partial\Omega} (n\times\Delta u^\alpha)\cdot (\beta(u^\alpha-v^\alpha)\times n)\\[3mm]
= & \displaystyle \int_\Omega (\nabla\times(\Delta u^\alpha))\cdot(\beta(u^\alpha-v^\alpha)\times n)-\int_\Omega \Delta u^\alpha\cdot \nabla\times(\beta(u^\alpha-v^\alpha)\times n)
\end{array}
\end{eqnarray}
where we have extended $\beta$ and $n$ smoothly to $\bar{\Omega}$. It follows from (4.7) that
$$||v^\alpha-u^\alpha||^2=(-\alpha\Delta u^\alpha, v^\alpha-u^\alpha)\leq \alpha||\Delta u^\alpha||\,||v^\alpha-u^\alpha|| $$
which yields
\begin{eqnarray}
||v^\alpha-u^\alpha||\leq \alpha||\Delta u^\alpha||
\end{eqnarray}
It follows from (6.57) and (6.68) that
\begin{eqnarray}
\left|\int_\Omega\Delta u^\alpha\cdot\nabla\times(\beta(u^\alpha-v^\alpha)\times n)\right| \leq \frac{1}{4}||\Delta u^\alpha||^2 + c||v^\alpha||^2_1
\end{eqnarray}
for suitably small $\alpha$. Using (6.68) again gives
\begin{eqnarray}
\begin{array}{rcl}
\displaystyle &&\left|\int_\Omega(\nabla\times\Delta u^\alpha)\cdot(\beta(u^\alpha-v^\alpha)\times n)\right| \\
 &\leq & \displaystyle \frac{1}{2}\alpha\int_\Omega |\nabla\times(\Delta u^\alpha)|^2 + \alpha^{-1}c||u^\alpha-v^\alpha||^2\\[3mm]
& \leq & \displaystyle\frac{1}{2}\alpha\int_\Omega |\nabla\times(\Delta u^\alpha)|^2 + c\alpha||\Delta u^\alpha||^2
\end{array}
\end{eqnarray}
Collecting (6.66),(6.67), (6.69) and (6.70) leads to
\begin{eqnarray}
||\Delta u^\alpha||^2 + \alpha||(\nabla\times)^3\,u^\alpha||^2 \leq c||v^\alpha||^2_2
\end{eqnarray}
for suitably small $\alpha$. This, together with the bound of $\partial_t\,u^\alpha$ in $H$,
implies the desired convergence in (6.50),(6.51). Thus the theorem is proved.

\subsection{Estimates on Convergence Rates}

Finally, we study the rates of convergence in the case of strong
solutions. We start with the case that the limiting Navier-Stokes
system has a strong solution.

\begin{thm} Let $v_0\in V$ and $v^0$ be the strong
solution to the Navier-Stokes equation with initial data $v_0$ on
any given finite interval $[0,T]$ with $T>0$. Then there exists a
$\alpha_0>0$ such that for each $\alpha\in(0,\alpha_0]$, the
LNS-$\alpha$ with the initial data $v_0$ has a unique strong
solution $(v^\alpha, u^\alpha)$ on the same interval $[0,T]$
satisfying
\begin{eqnarray}
\sup_{0\leq t\leq T} ||(v^\alpha,u^\alpha) - (v^0,v^0)||^2 + \int^T_0||(v^\alpha,u^\alpha) - (v^0,v^0)||^2_1(t)dt \leq c\alpha
\end{eqnarray}
\begin{eqnarray}
\sup_{0\leq t\leq T}||v^\alpha-v^0||^2_1 + \int^T_0 ||v^\alpha-v^0||^2_2 dt \leq c\alpha^{\frac{1}{2}}
\end{eqnarray}
with $c$ being a positive constant depending on $v^0$.
\end{thm}
\textbf{Proof:} Thanks to the local well-posedness of the strong
solution to the initial-boundary value problem for the
LNS-$\alpha$ and the standard continuation arguments, theorem 6.3.
will follow immediately from the following a priori estimates.

\begin{prop} Let $T_1\in(0,T]$ and $(u^\alpha,v^\alpha)$ be the strong solution
to the LANS-$\alpha$ with the initial data $v^0$ on the interval $[0,T_1]$ with the property that
\begin{eqnarray}
||v^\alpha||^2_1(t) +\int^t_0||v^\alpha||^2_2(\tau) d\tau \leq c_0 \qquad {\rm for} \ t\in[0,T_1]
\end{eqnarray}
with a positive constant $c_0$ depending only on $T$ and $v^0\in L^\infty(0,T:V)\cap L^\tau(0,T: W_\beta)$.
Then there exist uniform constants $\alpha_1$ and $c$ with the same dependence as $c_0$ such that
\begin{eqnarray}
\sup_{0\leq t\leq T_1}||(v^\alpha,u^\alpha)-(v^0,v^0)||^2+\int^{T_1}_0||(v^\alpha,u^\alpha)-(v^0,v^0)||^2_1(t) dt\leq c\alpha
\end{eqnarray}
\begin{eqnarray}
\sup_{0\leq t\leq T_1} ||(v^\alpha-v^0)||^2_1(t) + \int^{T_1}_0 ||v^\alpha-v^0||^2_2(t) dt \leq c\alpha^{\frac{1}{2}}
\end{eqnarray}
for $\alpha\in(0,\alpha_1]$.
\end{prop}

\begin{rem}
Assuming Proposition 6.4. for a moment, one can verify the a priori assumption
(6.74) by choosing
\[ c_0=1+4\sup_{0\leq t\leq T} ||v^0||^2_1 + \int^T_0||v^0||^2_2dt \]
and using (6.75) and (6.76) to choose $\alpha_0$. Thus (6.75),(6.76) hold for
$T_1=T$. This yields theorem 6.3. immediately. It remains to verify proposition 6.4.
\end{rem}
\textbf{Proof of Proposition 6.4.}
Set $w = v^{\alpha} - v^0$. Then it holds that for $t\in[0,T_1]$,
\begin{equation}
\frac{d}{dt}\|w\|^2 (t) + 2a_\beta(w,w)(t) + 2 (B (v^{\alpha},u^\alpha) -
B(v^0,u^0), w)(t) = 0
    \end{equation}
Note that
$$2(B(v^\alpha,u^\alpha)-B(v^0,v^0),w) = 2(B(v^\alpha,v^\alpha)-B(v^0,v^0),w) + 2\alpha((\nabla\times v^\alpha)\times P(\Delta u^\alpha),w)$$
$$|(B(v^\alpha,v^\alpha)-B(v^0,v^0),w)| = |(B(w,v^0),w)| \leq c||v^0||_2 \, ||w||_1 \, ||w||$$
$$|(\nabla\times v^\alpha)\times P(\Delta u^\alpha),w|\leq c||P(\Delta u^\alpha)|| \, ||v^\alpha||_2 \, ||w||_1$$
It follows that for all $t\in[0,T_1]$
\begin{equation}
\frac{d}{dt}\|w\|^2(t) + a_\beta(w,w)(t) \leq c(1+\|v^0\|^2_2) \|w\|^2 +
c\alpha^2\|P(\Delta u^{\alpha})\|^2 \, \|v^{\alpha}\|_2^2
\end{equation}
Due to (6.59), (6.60) and (6.74), one has
\begin{eqnarray}
\alpha||P(\Delta u^\alpha)||^2 \leq c c_0 \qquad {\rm for\ all} \ t\in[0,T_1]
\end{eqnarray}
which, together with (6.78), yields
\begin{equation}
\frac{d}{dt}\|w\|^2(t) + a_\beta(w,w)(t) \leq c(1+\|v^0\|^2_2)\,\|w\|^2(t) +
c c_0\alpha\|v^{\alpha}\|_2^2 (t)
\end{equation}
Since $w(0)=0$, so Gronwall's inequality leads to
\begin{eqnarray}
\begin{array}{rcl}
\displaystyle &&\sup_{0\leq t\leq T_1} ||w(t)||^2 + \int^{T_1}_0 ||w||^2_1(t) dt\\
& \leq & \displaystyle (c\alpha c_0\,e^{c\int^{T_1}_0(1+||v^0||^2_2(t))dt}) \, \int^{T_1}_0 ||v^\alpha||^2_2(t) dt\\[3mm]
& \leq & c\alpha\,c^2_0 e^{cc_1} \ \equiv \ c_2\alpha
\end{array}
\end{eqnarray}
with $c_1$ depending only $T$ and the $L^2(0,T:H^2)$-norm of
$v^0$.\\

To prove (6.75) for $u^\alpha$, we note that (4.3) implies that
\begin{eqnarray}
||u^\alpha-v^0||^2 + \alpha a_\beta(u^\alpha,u^\alpha-v^0)=(v^\alpha-v^0,u^\alpha-v^0)
\end{eqnarray}
\begin{eqnarray}
a_\beta(u^\alpha-v^0,u^\alpha-v^0)+\alpha(P(\Delta u^\alpha),\Delta(u^\alpha-v^0)) = a_\beta (v^\alpha-v^0,u^\alpha-v^0)
\end{eqnarray}
It follows from (6.82), (6.83), and (6.74) that
\begin{eqnarray}
\begin{array}{rcl}
\displaystyle ||u^\alpha-v^0||^2 + a_\beta(u^\alpha-v^0,u^\alpha-v^0) & \leq & \displaystyle ||v^\alpha-v^0||^2 + \alpha a_\beta (v^0,v^0) \\
& \leq & c_2\alpha +cc_0\alpha
\end{array}
\end{eqnarray}
Since
\begin{eqnarray*}
\begin{array}{rcl}
\alpha(P(\Delta u^\alpha),\Delta(u^\alpha-u^0)) & = & \alpha||P(\Delta u^\alpha)||^2 - \alpha(P(\Delta u^\alpha),\Delta v^0)\\
& \geq & \displaystyle \frac{1}{2} \alpha||P(\Delta u^\alpha)||^2 - \frac{\alpha}{2}||\Delta v^0||^2
\end{array}
\end{eqnarray*}
This, together with (6.84), shows that
\begin{eqnarray}
a_\beta(u^\alpha-v^0, u^\alpha-v^0)+\alpha||P(\Delta u^\alpha)||^2 \leq a_\beta (v^\alpha-v^0,v^\alpha-v^0)+\alpha||\Delta v^0||^2
\end{eqnarray}
Hence, one obtains from (6.82), (6.84), and (6.85) that
$$\sup_{0\leq t\leq T_1} ||u^\alpha-v^0||^2 + \int^{T_1}_0||u^\alpha-v^0||^2_1\,dt \leq c_3\,\alpha$$
This completes the verification of (6.75). It remains to prove
(6.76).
\\

By the definition of strong solutions, one has that for a.e. $t\in[0,T_1]$,
\begin{eqnarray}
\frac{d}{dt} a_\beta (w,w) + 2||P\Delta w||^2 + 2(B(v^\alpha,u^\alpha)-B(v^0,v^0),P\Delta w)=0
\end{eqnarray}
Rewrite the last term on the left hand side above as
$$2(B(v^\alpha,u^\alpha)-B(v^0,v^0),-P(\Delta w)) = I+II$$
One can estimate each term as
$$|I_1|=|2(B(v^\alpha,v^\alpha)-B(v^0,v^0), -P(\Delta w))|\leq c(||v^\alpha||_2 + ||v||_2)||w||_1\, ||P\Delta w||$$
$$|II_2|=|2\alpha((\nabla\times v^\alpha)\times(P\Delta u^\alpha),P\Delta w)|\leq c^\alpha||P\Delta u^\alpha||_1 \, ||v^\alpha||^{\frac{1}{2}}_1 \, ||v^\alpha||^{\frac{1}{2}}_2 \, ||P\Delta w||$$
It follows that for a.e. $t\in[0,T_1]$,
\begin{eqnarray}
\begin{array}{rcl}
& & \displaystyle \frac{d}{dt} \tilde{a}_\beta (w,w)(t) +||P\Delta w||^2(t)\\[1mm]
& \leq & c(1+||v^0||^2_2 + ||v^\alpha||^2_2) \tilde{a}_\beta(w,w)\\[1mm]
& & + c\alpha^2 ||P\Delta u^\alpha||^2_1 \, ||v^\alpha||_1 \, ||v^\alpha||_2 + c\alpha^2 ||P\Delta u^\alpha||^2 \, ||v^\alpha||^2_2
\end{array}
\end{eqnarray}
As a consequence of (6.57), (6.74) and
$$\alpha||P(\Delta u^\alpha)||_1 = ||u^\alpha - v^\alpha||_1$$
one has
\begin{eqnarray}
\begin{array}{cl}
& c\alpha^2 \, ||P\Delta u^\alpha||_1 \, ||v^\alpha||_1 \, ||v^\alpha||_2\\
\leq & c||v^\alpha-u^\alpha||^2_1 \, ||v^\alpha||_1 \, ||v^\alpha||_2\\
\leq & c||v^\alpha||_1 \, ||v||^\alpha_2 \, \tilde{a}_\beta(w,w) + c||u^\alpha - v^0||^2_1 \, ||v^\alpha||_1 \, ||v^\alpha||_2
\end{array}
\end{eqnarray}
It follows from (6.87), (6.88), and (6.79) that
\begin{eqnarray}
\begin{array}{rcl}
& & \displaystyle \frac{d}{dt} \tilde{a}_\beta (w,w) +\nu||A_\beta w||^2\\[1mm]
& \leq & c(1+||v^0||^2_2 + ||v^\alpha||^2_2) \tilde{a}_\beta(w,w)\\[1mm]
& & + cc^{\frac{1}{2}}_0 \, ||u^\alpha-v^0||^2_1 \, ||v^\alpha||_2 + cc_0\alpha ||v^\alpha||^2_2
\end{array}
\end{eqnarray}
Consequently, we can get
\begin{eqnarray*}
\begin{array}{cl}
& \displaystyle \sup_{0\leq t\leq T_1} ||w(t)||^2_1 + \int^{T_1}_0 ||w(t)||^2_2 dt\\[1mm]
\leq & \displaystyle \int^{T_1}_0 e^{c\int^{T_1}_0(1+||v_0||^2_2+||v^\alpha||^2_2)dt} \ [ cc_0\alpha||v^\alpha(t)||^2_2 + cc^{\frac{1}{2}}_0||u^\alpha-v^0||^2_1 \, ||v^\alpha||_2] \,dt\\
\leq & \displaystyle c_1\,\alpha+c_1\int^{T_1}_0 ||u^\alpha-v^0||^2_1 \, ||v^\alpha||_2 \,dt\\
\leq & \displaystyle c_1\,\alpha+c_1 \left(\int^{T_1}_0 ||u^\alpha-v^0||^4_1 dt \right)^{\frac{1}{2}} \left(\int^{T_1}_0 ||v^\alpha||^2_2 dt \right)^{\frac{1}{2}}\\
\leq & \displaystyle c_1\,\alpha+c_1 (||v^\alpha||_1+||v^0||^2_1) \left(\int^{T_1}_0||u^\alpha-v^0||^2_1 dt \right)^{\frac{1}{2}}
\leq  c_1\,\alpha^{\frac{1}{2}}
\end{array}
\end{eqnarray*}
where we have used (6.75). Thus (6.76) holds, and the proposition
is proved.
\begin{rem} It is not clear to us whether the stronger estimate as (6.73)
holds for $u^\alpha$ under the assumptions in theorem 6.3..
However, under the additional assumption that the strong solution $v^0\in L^\infty([0,T],H^2)$, there holds also
\begin{eqnarray}
\sup_{0\leq t\leq T} ||u^\alpha-v^0||^2_1 + \int^T_0||u^\alpha-v^0||^2_2 \, dt\leq c\alpha^{\frac{1}{2}}
\end{eqnarray}
This follows from
\begin{eqnarray*}
\begin{array}{cl}
& a_\beta(u^\alpha-v^0,u^\alpha-v^0) + \alpha(P\Delta(u^\alpha-v^0),P\Delta(u^\alpha-v^0))\\
= & a_\beta(v^\alpha-v^0,u^\alpha-v^0)+\alpha(-P\Delta v^0,P\Delta(u^\alpha-v^0))
\end{array}
\end{eqnarray*}
(due to (6.83)) and (6.73).
\end{rem}

\section{Concluding Remarks}

We conclude this paper with a few remarks on related issues.

\begin{rem} In exact same way, we can study the boundary value problem
of LNS-$\alpha$ with NSB:
\begin{eqnarray}
&& \partial_t v - \Delta v + \nabla\times v \times u +
\nabla p=0 \ {\rm in} \ \Omega\\
&& \nabla\cdot v = 0 \ {\rm in} \ \Omega\\
&& u - \alpha \Delta u +
\nabla \tilde{p} = v \ {\rm in} \ \Omega\label{a1}\\
&& \nabla\cdot u = 0 \ {\rm in} \ \Omega\\
&& v \cdot n = 0,\  2(S(v)n)_\tau = -\gamma v_\tau \ {\rm on} \ \partial\Omega\\
&& u \cdot n = 0,\  2(S(u)n)_\tau = -\gamma u_\tau \ {\rm on} \
\partial\Omega
\end{eqnarray}
The functional setting is similar to that of
(4.1)-(4.6), and all the results stated in section 4-6 are also valid.
\end{rem}
\begin{rem} The non-homogenous boundary value problems
of LNS-$\alpha$ with VSB:
\begin{eqnarray}\label{7.1}
&& \partial_t v - \Delta v + \nabla\times v \times u +
\nabla p=0 \ {\rm in} \ \Omega\\
&& \nabla\cdot v = 0 \ {\rm in} \ \Omega\\
&& u - \alpha \Delta u +
\nabla \tilde{p} = v \ {\rm in} \ \Omega\label{a1}\\
&& \nabla\cdot u = 0 \ {\rm in} \ \Omega\\
&& v \cdot n = 0,\  n\times \nabla\times v = \beta v + b \ {\rm
on} \ \partial\Omega\\\label{7.2}
&& u \cdot n = 0,\  n\times
\nabla\times u = \beta v + b \ {\rm on} \ \partial\Omega
\end{eqnarray}
can also be considered by using a homogenous method to
reduce it into
\begin{eqnarray}
&& \partial_t v - \Delta v + \nabla\times v \times u +
\nabla p = \xi \ {\rm in} \ \Omega\\
&& \nabla\cdot v = 0 \ {\rm in} \ \Omega\\
&& u - \alpha \Delta u +
\nabla \tilde{p} = v + \eta_\alpha \ {\rm in} \ \Omega\\
&& \nabla\cdot u = 0 \ {\rm in} \ \Omega\\
&& v \cdot n = 0,\  n\times \nabla\times v = \beta v \ {\rm on} \ \partial\Omega\\
&& u \cdot n = 0,\  n\times \nabla\times u = \beta u \ {\rm on} \
\partial\Omega
\end{eqnarray}
for some $\xi$ and $\eta_\alpha$ as was done for the steady homogenous case in Section 3.
Similarly, the non-homogenous boundary
value problems for LNS-$\alpha$ with NSB may be established too.
\end{rem}
\begin{rem} In the functional settings, the parameters associated with the velocity $v$ and the filter $u$ can be different, and
different type boundary conditions, VSB or NSB, may be also allowed. However, in this case the analysis in the global existence and the
vanishing $\alpha$ limit seems very difficult since (\ref{5.2}) does not hold, there
are some boundary terms arising, and the energy estimate in (\ref{5.3}) depends on $\alpha$.
Yet the local well-posedness theory can be established by the method discussed in this paper.
\end{rem}
\begin{rem} Our approaches works also for other $\alpha$ models. For instance, one can consider the following Leray $\alpha$
model:
\begin{eqnarray}
&& \partial_t v - \Delta v + u\cdot \nabla v  +
\nabla p=0 \ {\rm in} \ \Omega\\
&& \nabla\cdot v = 0 \ {\rm in} \ \Omega\\
&& u - \alpha \Delta u +
\nabla \tilde{p} = v \ {\rm in} \ \Omega\\
&& \nabla\cdot u = 0 \ {\rm in} \ \Omega\\
&& v \cdot n = 0,\ n\times (\nabla\times v) = \beta v  \ {\rm on} \ \partial\Omega\\
&& u \cdot n = 0,\  2(S(u)n)_\tau = -\gamma u_\tau \ {\rm on} \
\partial\Omega
\end{eqnarray}
which allowed different boundary conditions between the velocity $v$
and the filter $u$. In fact, this model is easier to analyze than the LNS-$\alpha$
since it has the following energy equation
\begin{equation}
  \frac{d}{dt}\|v\|^{2} + 2 a_{\beta}(v,v)  = 0
    \end{equation}
which yields the global existence directly, and the corresponding
convergence result is better both in $v^\alpha$ and $u^\alpha$ than that in
theorem 5.1..
\end{rem}

\end{document}